\documentclass[final]{etna_arxiv}
\usepackage{algorithmicx}
\usepackage{algorithm}
\usepackage{algpseudocode}
\usepackage{amsmath,amssymb}
\usepackage{mathtools}
\usepackage[dvipsnames]{xcolor}
\usepackage{soul} 
\usepackage{adjustbox}
\usepackage{graphicx}
\usepackage{subcaption}

\renewcommand{\algorithmicrequire}{\textbf{Input: }}
\renewcommand{\algorithmicensure}{\textbf{Output: }}

\newcommand{\Atilde}{\widetilde{A}}
\newcommand{\Ahat}{\widehat{A}}
\newcommand{\btilde}{\tilde{b}}
\newcommand{\Cbar}{\bar{C}}
\newcommand{\ctilde}{\tilde{c}}

\newcommand{\Hbar}{\bar{H}}
\newcommand{\Lbar}{\bar{L}}
\newcommand{\Rbar}{\bar{R}}
\newcommand{\rbar}{\bar{r}}
\newcommand{\sbar}{\bar{s}}
\newcommand{\rtilde}{\tilde{r}}
\newcommand{\Ubar}{\bar{U}}
\newcommand{\Vhat}{\widehat{V}}
\newcommand{\Vbar}{\bar{V}}

\newcommand{\Zbar}{\bar{Z}}

\newcommand{\Zhat}{\widehat{Z}}
\newcommand{\wbar}{\bar{w}}
\newcommand{\ybar}{\bar{y}}
\newcommand{\xbar}{\bar{x}}

\newcommand{\betabar}{\bar{\beta}}

\setbibdata{1}{xx}{46}{2017} 

\hypersetup{%
    pdftitle={Document title},
    pdfauthor={John Doe, Erwin Schr\"odinger},
    pdfkeywords={a keyword, another keyword, and another one}
    }

\title{The stability of split-preconditioned FGMRES in four precisions\thanks{%
 Co-funded by the Exascale Computing Project (17-SC-20-SC), a collaborative
effort of the U.S. Department of Energy Office of Science and the National Nuclear Security Administration, and by the European Union (ERC, inEXASCALE, 101075632). Views and opinions expressed are however those of the authors only and do not necessarily reflect those of the European Union or the European Research Council. Neither the European Union nor the granting authority can be held responsible for them. 
}}

\author{Erin Carson \footnotemark[2]
        \and Ieva Dau\v{z}ickait\.{e}\footnotemark[2]}

\shorttitle{Mixed precision FGMRES} 
\shortauthor{E.~Carson AND I.~Dau\v{z}ickait\.{e}}

\begin{document}

\maketitle

\renewcommand{\thefootnote}{\fnsymbol{footnote}}

\footnotetext[2]{Faculty of Mathematics and Physics, Charles University, Sokolovsk\'{a} 366/84, 180 00, Prague, Czech Republic. email: carson@karlin.mff.cuni.cz, dauzickaite@karlin.mff.cuni.cz.}

\begin{abstract}
We consider the split-preconditioned FGMRES method in a mixed precision framework, in which four potentially different precisions can be used for computations with the coefficient matrix, application of the left preconditioner, application of the right preconditioner, and the working precision. Our analysis is  applicable to general preconditioners. We obtain bounds on the backward and forward errors in split-preconditioned FGMRES. Our analysis further provides insight into how the various precisions should be chosen; under certain assumptions, a suitable selection guarantees a backward error on the order of the working precision. 
\end{abstract}



\section{Introduction}
We consider the problem of solving a linear system of equations
\begin{equation}\label{eq:Ax=b}
    Ax=b,
\end{equation}
where $A \in \mathbb{R}^{n \times n}$ is nonsymmetric and $x,b \in \mathbb{R}^n$. When $A$ is large and sparse, the iterative generalised minimal residual method (GMRES) or its flexible variant (FGMRES) are often used for solving \eqref{eq:Ax=b}; see, for example, \cite{saad2003iterative}. 
In these and other Krylov subspace methods, preconditioning is an essential ingredient. Given a preconditioner $P = M_L M_R$, the problem \eqref{eq:Ax=b} is transformed to
\begin{align}
   M_L^{-1} A M_R^{-1} \tilde{x} = \, &  M_L^{-1} b, \label{eq:Ax=b_split_precond}\\
   \textrm{where } M_R^{-1} \tilde{x} = \, & x. \nonumber
\end{align}
Note that a particular strength of FGMRES is that it allows the right preconditioner to change throughout the iterations. Although for simplicity, we consider the case here where the preconditioners are static, our results could be extended to allow dynamic preconditioning.

The emergence of mixed precision hardware has motivated work in developing mixed precision algorithms for matrix computations; see, e.g., the recent surveys \cite{abdelfattah2021survey,higham_mary_2022}. Modern GPUs offer double, single, half, and even quarter precision, along with specialized tensor core instructions; see, e.g., \cite{h100}. The use of lower precision can offer significant performance improvements, although this comes at a numerical cost. With fewer bits, we have a greater unit roundoff and a smaller range of representable numbers. The goal is thus to selectively use low precision in algorithms such that performance is improved without adversely affecting the desired numerical properties.

Mixed precision variants of GMRES and FGMRES with different preconditioners have been proposed and analyzed in multiple papers. 
Arioli and Duff \cite{arioli2009using} analyzed a two-precision variant of FGMRES in which the right-preconditioner is constructed using an LU decomposition computed in single precision and applied in either single or double precision, and other computations are performed in double precision. They proved that in this setting, a backward error on the order of double precision in attainable. 
The authors of \cite{lindquist2020improving} develop a mixed precision variant of left-preconditioned GMRES in a mix of single and double precisions, requiring only a few operations to be performed in double precision. Their numerical experiments show that they can obtain a backward error to the level of double precision. 
Variants of left-preconditioned GMRES using various numbers of precisions have been analyzed as inner solvers within GMRES-based iterative refinement for solving linear systems of equations;  
Vieubl\'{e} \cite{vieuble_thesis} analyzed left-preconditioned GMRES in four precisions with a general preconditioner, following the earlier works \cite{carson2017new} and \cite{h:21} which analyzed left-preconditioned GMRES with an LU preconditioner in two and three precisions, respectively. In general, different precisions can be used for computing the preconditioner, matrix-vector products with $A$, matrix-vector products or solves with the general preconditioner(s), and the remaining computations. 
We refer the readers to the recent surveys \cite{abdelfattah2021survey,higham_mary_2022} for other examples.

The structure of some problems and/or application requirements makes it desirable to construct and apply a split-pre\-cond\-i\-tion\-er rather than left or right ones alone. For example, the condition number of a split-preconditioned matrix can be significantly smaller than when the same preconditioner is applied on the left \cite[Section 3.2]{carson2020three}. Such preconditioning is usually used for symmetric problems solved via short recurrence symmetric solvers such as MINRES. However, MINRES may be inferior compared to GMRES when a high accuracy solution for an ill-conditioned problem is required \cite{carson2020three}. We also emphasize that analyzing split-preconditioning provides a uniform framework for analyzing cases with full left- or full right-preconditioning.
The stability of split-preconditioned GMRES and FGMRES has not been analyzed in either uniform or mixed precision. 
The work \cite{arioli2007note} showed that uniform precision FGMRES with a specific right-preconditioner is backward stable while this is not the case for GMRES and that FGMRES is more robust than GMRES. We thus focus on split-preconditioned FGMRES in this paper and develop a mixed precision framework allowing for four potentially different precisions for the following operations:
computing matrix-vector products with $A$, applying the left-preconditioner $M_L$, applying the right-preconditioner $M_R$, and all other computations. FGMRES computes a series of approximate solutions $x_k$ from Krylov subspaces to \eqref{eq:Ax=b}. The Arnoldi method is employed to generate the basis for the Krylov subspaces like in GMRES, but FGMRES stores the right-preconditioned basis as well. The particular algorithm is shown in Algorithm~\ref{alg:fgmres}. Our analysis considers general preconditioners, only requiring an assumption on the error in applying its inverse to a vector, and is thus widely applicable. 

The paper is outlined as follows. We bound the backward errors in Section~\ref{sec:finite_prec_analysis} while also providing guidance for setting the four precisions such that backward error to the desired level is attainable. To make the results of the analysis more concrete, in Section~\ref{sec:example_LU_precond} we bound the quantities involved for the example of LU preconditioners, and then present a set of numerical experiments on both dense problems and problems from SuiteSparse \cite{SuiteSparse}. In Section~\ref{sec:conclusions} we make concluding remarks.

\begin{algorithm}
\caption{Split-preconditioned FGMRES for solving \eqref{eq:Ax=b} in four precisions}\label{alg:fgmres}
\algorithmicrequire  matrix $A \in \mathbb{R}^{n \times n}$, right hand side $b \in \mathbb{R}^n$, preconditioner $P = M_L M_R$,  maximum number of iterations \textit{maxit}, convergence tolerance $\tau$, precisions $u$, $u_A$, $u_L$ and $u_R$\\
\algorithmicensure approximate solution $x_k$
\begin{algorithmic}
\State  \textrm{initialize } $x_0$ 
\State $t = Ax_0$ \Comment{$u_A$}
\State $t^{(p)} = M_L^{-1} t$ \Comment{$u_L$}
\State $ b^{(p)} = M_L^{-1}b$ \Comment{$u_L$}
\State $r_0 = b^{(p)} - t^{(p)}$ \Comment{$u$}
\State $\beta = \Vert r_0 \Vert$; $v_1 = r_0 / \beta$ \Comment{$u$}  
\State  $k=0$; \textit{convergence} = \textit{false}
\While {\textit{convergence} = \textit{false} and $k < maxit$}
\State $k = k +1$ 
\State $z_k = M_R^{-1} v_k$ \Comment{$u_R$}
\State $s = A z_k$ \Comment{$u_A$}
\State $w = M_L^{-1} s$ \Comment{$u_L$}
\For {$i = 1,\dots,k$}
\State $h_{i,k} = v_i^T w$ \Comment{$u$} 
\State $w = w - h_{i,k} v_i$ \Comment{$u$} 
\EndFor
\State $h_{k+1,k} = \Vert w \Vert$ \Comment{$u$} 
\State $Z_k = [z_1, \dots, z_k]$; $H_k = \{h_{i,j}\}_{1 \leq i \leq j+1; 1 \leq j \leq k}$
\State $y_k = \textrm{arg} \min_y \Vert \beta e_1 - H_k y \Vert$ \Comment{$u$} 
\If {$\Vert \beta e_1 - H_k y_k \Vert \leq \tau  \beta $}
\State $x_k = x_0 + Z_k y_k $  \Comment{$u$ } 
\State $t = Ax_k$  \Comment{$u_A$}  
\State $r = b - t$  \Comment{$u$ }  
\State \textit{convergence} = \textit{true}
\Else
\State $v_{k+1} = w/h_{k+1,k}$; $V_{k+1} = [v_1, \dots, v_{k+1}]$  \Comment{$u$} 
\EndIf
\EndWhile
\end{algorithmic}
\end{algorithm}

\section{Finite precision analysis of FGMRES in four precisions}\label{sec:finite_prec_analysis}
From the Rigal-Gaches Theorem (see \cite[Theorem 7.1]{high:ASNA2}), the normwise relative backward error is given by 
\[
\min\{\varepsilon:(A+\Delta A)x_k = b+\Delta b, \|\Delta A\|\leq \varepsilon \|A\|, \|\Delta b\|\leq \varepsilon\|b\| \} = \frac{\| r_k\|}{\|A\|\|x_k\|+\|b\|},
\]
where $r_k=b-Ax_k$. We aim to bound this quantity when $x_k$ is the approximate solution produced by Algorithm~\ref{alg:fgmres}. To account for various ways in which the preconditioner can be computed and some constraints on $A$ resulting in the need for different precisions, we assume that
\begin{itemize}
    \item computations with $A$ are performed in precision with unit roundoff $u_A$;
    \item computations with $M_L$ are performed in precision with unit roundoff $u_L$;
    \item computations with $M_R$ are performed in precision with unit roundoff $u_R$;
    \item the precision for other computations (the working precision) has unit roundoff $u$.
\end{itemize}
Note that when these precisions differ, some conversion between them is required. This may be done implicitly or explicitly depending on the particular precisions and underlying hardware and software. We also highlight that in Algorithm~\ref{alg:fgmres} computing $r_0$ requires computing $M_L^{-1} A x_0$ and $M_L^{-1} b$ separately instead of the usual $M_L^{-1}(b-Ax_0)$; the prior choice is more computationally expensive and we only do this to avoid terms $u \Vert \vert M_L^{-1}\vert \vert M_L \vert \Vert$ in the analysis. If $x_0=0$, then $M_L^{-1}$ needs to be applied only once.

Using the approach in \cite{vieuble_thesis}, we assume that the application of $M_L^{-1}$ and $M_R^{-1}$ can be computed in a way such that
\begin{align}
    fl(M_L^{-1} w_j) = & M_L^{-1} w_j + \Delta M_{L,j} w_j, \quad  \lvert \Delta M_{L,j} \rvert \leq c(n) u_L E_{L,j}, \label{eq:Ml_assump} \\
    fl(M_R^{-1} w_j) = & M_R^{-1} w_j + \Delta M_{R,j} w_j, \quad  \lvert \Delta M_{R,j} \rvert \leq c(n) u_R E_{R,j}, \label{eq:Mr_assump}
\end{align}
where $fl(\cdot)$ denotes the quantity computed in floating point arithmetic, $E_{L,j}$ and $E_{R,j}$ have positive entries, $w_j \in \mathbb{R}^n$, and $c(n)$ is a constant that depends on $n$ only. We define
\begin{equation*}
    \Atilde \coloneqq M_L^{-1} A \quad \text{ and } \quad \btilde \coloneqq M_L^{-1} b
\end{equation*}
and assume that matrix-vector products with $\Atilde$ can be computed so that
\begin{equation*}
    fl(\Atilde z_j) = (M_L^{-1} +  \Delta M_{L,j}) (A +  \Delta A_j) z_j.
\end{equation*}
Denoting
\begin{equation*}
    u_A \psi_{A,j} = \frac{\Vert M_L^{-1} \Delta A_j z_j \Vert}{\Vert \Atilde \Vert \Vert z_j \Vert} \quad \text{ and } \quad   u_L \psi_{L,j} = \frac{\Vert \Delta M_{L,j} A z_j\Vert}{\Vert \Atilde \Vert \Vert z_j \Vert},
\end{equation*}
where here and in the rest of the paper $\Vert \cdot \Vert$ denotes the 2-norm, and ignoring the second order terms, we can write 
\begin{gather*}
    fl(\Atilde z_j) \approx \Atilde z_j + f_j,\\
   \text{where } \Vert f_j \Vert \leq (u_A \psi_{A,j} + u_L \psi_{L,j})\Vert \Atilde \Vert \Vert z_j \Vert.
\end{gather*}

In the following, a standard error analysis approach is used, e.g., \cite{high:ASNA2}, and  we closely follow the analysis in \cite{arioli2007note} and \cite{arioli2009using}. The analysis is performed in the following stages:
\begin{enumerate}
    \item Bounding the computed quantities in the modified Gram-Schmidt (MGS) algorithm that returns $C^{(k)} = \begin{bmatrix} \btilde - \Atilde x_0 &  \Atilde Z_k \end{bmatrix} = V_{k+1} R_k$, where $ V_k^T V_k =   I_k$, $R_k = \, \begin{bmatrix} \beta e_1 & H_k \end{bmatrix}$ and $e_1 = \begin{bmatrix} 1 & 0 & \dots & 0\end{bmatrix}^T$.
    \item Solving the least-squares problem 
    \begin{equation}\label{eq:LS_problem}
       y_k = \textrm{arg} \min_y \Vert \beta e_1 - H_k y \Vert
    \end{equation}
    via QR employing Givens rotations and analyzing its residual.
    \item Computing $x_k = x_0 + Z_k y_k$.
    \item Bounding $\Vert y_k \Vert$.
\end{enumerate}
Throughout the paper computed quantities are denoted with bars, that is, $\Cbar^{(k)}$ is the computed $C^{(k)}$, 
and $\kappa(A) = \Vert A \Vert \Vert A^{\dagger} \Vert$ is the 2-norm condition number of $A$. The second order terms in $u_A$, $u_L$, $u_R$, and $u$ are ignored. We drop the subscripts $j$ for $E_{L,j}$, $E_{R,j}$, $\Delta M_{L,j}$, $\Delta M_{R,j}$, $ \Delta A_j $, $\psi_{A,j}$ and $\psi_{L,j}$ and replace these quantities by their maxima over all $j$. 
It is assumed that no overflow or underflow occurs. We present the main result here and refer the reader to Appendix~\ref{app:backward_err_proof} for the proof.

\begin{theorem}\label{th:backward_error}
Let $\xbar_k$ be the approximate solution to \eqref{eq:Ax=b_split_precond} computed by Algorithm~\ref{alg:fgmres}. Under the assumptions \eqref{eq:Ml_assump}, \eqref{eq:Mr_assump}, 
\begin{gather}
    2.12(n+1)u < 0.01 \quad \text{and} \quad c_0(n) u \kappa(C^{(k)}) < 0.1 \ \forall k, \label{eq:assumptions_nu_ukappaC} \\
     \lvert \sbar_k \rvert < 1 - u \quad \forall k, \label{eq:assumption_sk}
\end{gather}
where $c_0(n) = 18.53 n^{3/2}$ and $\sbar_k$ are the sines computed for the Givens rotations, and
\begin{equation}\label{eq:rho_def_assump}
    \rho \coloneqq 1.3 c_{13}(n,k) \Vert M_R \Vert \left( u \Vert \Zbar_k \Vert + u_R \Vert E_R \Vert \right) < 1,
\end{equation}
where $c_{13}(n,k)$ is defined in  \eqref{eq:yk_norm}, the residual for the left-preconditioned system is bounded by
\begin{equation}\label{eq:residual_bound}
\Vert \btilde - \Atilde  \xbar_k  \Vert \lesssim   \frac{1.3 c(n,k)}{1 - \rho} \left( \zeta_1 + \zeta_2 \right),
\end{equation}
where
\begin{gather}
    \zeta_1 \coloneqq \left(u + u_L \Vert E_L M_L \Vert \right) \Vert \btilde \Vert , \quad\text{and} \nonumber \\
    \zeta_2 \coloneqq \left( u + u_A \psi_A + u_L \psi_L \right) \Vert \Atilde \Vert 
 \left( \Vert \Zbar_k \Vert \Vert M_R(\xbar_k - \xbar_0) \Vert + \Vert \xbar_0 \Vert \right) \label{eq:zeta2_def}
\end{gather}
and the normwise relative backward error for the left-preconditioned system 
is bounded by
\begin{equation}\label{eq:backward_error}
   \frac{ \Vert \btilde - \Atilde  \xbar_k  \Vert}{\Vert \btilde \Vert + \Vert \Atilde \Vert \Vert \xbar_k  \Vert} \lesssim   \frac{1.3 c(n,k)}{1 - \rho} \zeta,
\end{equation}
where
\begin{equation}\label{eq:zeta}
    \zeta \coloneqq \frac{\zeta_1 + \zeta_2 }{\Vert \btilde \Vert + \Vert \Atilde \Vert \Vert \xbar_k  \Vert }.
\end{equation}

\end{theorem}
We expect \eqref{eq:backward_error} to be dominated by $\zeta_2$, mainly due to the term $\Vert \Zbar_k \Vert \Vert M_R(\xbar_k - \xbar_0) \Vert$. As observed in \cite{arioli2007note,arioli2009using} and in our experiments (Sections~\ref{sec:numerics_dense} - \ref{sec:numerics_sparse}),  $\Vert \Zbar_k \Vert$ remains small in early iterations, but can be large if many iterations are needed for convergence. We expect the quantity $\Vert M_R(\xbar_k - \xbar_0) \Vert$ to aid in partially mitigating the size of $\Vert \Zbar_k \Vert$, so that $\zeta_2$ still gives good guarantees for the backward error. Note that if we were to obtain $\Vert \Atilde \Vert \Vert \xbar_k  \Vert $ in $\zeta_2$ by using $\Vert M_R(\xbar_k - \xbar_0) \Vert \leq \Vert M_R \Vert (\Vert \xbar_k \Vert + \Vert \xbar_0 \Vert)$, then we would introduce the term $\Vert M_R \Vert \Vert \Zbar_k \Vert$. Depending on the preconditioner, $\Vert M_R \Vert$ can be close to $\Vert A \Vert$ and for some problems $\Vert M_R \Vert \Vert \Zbar_k \Vert$ can grow rapidly, thus making \eqref{eq:backward_error} a large overestimate. We comment on how \eqref{eq:residual_bound} compares with other bounds for FGMRES available in the literature in the following section. The condition \eqref{eq:rho_def_assump}, the quantities $\psi_A$, $\psi_L$, and the role of different precisions are discussed in Section~\ref{sec:precisions}.

Bound \eqref{eq:backward_error} can be formulated with respect to the the original system, that is, without the left preconditioner, using inequalities $\Vert b - A  \xbar_k  \Vert \leq \Vert M_L \Vert \Vert \btilde - \Atilde  \xbar_k  \Vert$ and $\Vert b \Vert + \Vert A \Vert \Vert \xbar_k  \Vert \geq  (\Vert \btilde \Vert + \Vert \Atilde \Vert \Vert \xbar_k  \Vert )/\Vert M_L^{-1} \Vert$. Alternatively, we can use the fact that the relative backward error is bounded by the relative forward error (see Section~\ref{sec:forward_error}). We state the bound in the following corollary.   
\begin{corollary}\label{cor:be_original_system}
If the conditions in Theorem~\ref{th:backward_error} are satisfied, then the normwise relative backward error for the system \eqref{eq:Ax=b} is bounded by
\begin{equation*}
    \frac{ \Vert b - A  \xbar_k  \Vert}{\Vert b \Vert + \Vert A \Vert \Vert \xbar_k  \Vert} \lesssim   \frac{1.3 c(n,k)}{1 - \rho} \zeta  \min\{ \kappa(M_L), \kappa(\Atilde)\}.
\end{equation*}
\end{corollary}
The condition number of the left preconditioner weakens the result, yet for some preconditioners $\kappa(M_L)$ can be expected to be small, for example when $LU$ decomposition is used and $M_L=L$. Note that a small backward error with respect to the preconditioned system and small $\kappa(\Atilde)$ implies small backward error with respect to the original system.

\subsection{Comparison with existing bounds}
We wish to compare our result with bound (5.6) in \cite{arioli2007note} for FGMRES with a general right-preconditioner and bound (3.32) in \cite{arioli2009using} for FGMRES right-preconditioned with $LU$ factorization computed in single precision. We set $M_L = I$, $u = u_A$, then $u_L = 0$ and $u_A \psi_A = u$. The bound \eqref{eq:residual_bound} becomes
\begin{equation*}
   \Vert b - A  \xbar_k  \Vert \lesssim \frac{1.3 c(n,k)u }{1 - \rho}  \left( \Vert b \Vert + \Vert A \Vert  \left( \Vert \Zbar_k \Vert \Vert M_R(\xbar_k - \xbar_0) \Vert + \Vert \xbar_0 \Vert \right) \right).
\end{equation*}
We thus recover the bound (5.6) in \cite{arioli2007note}, but ignoring the term $u^2 \Vert \xbar_0 \Vert$ and with a slightly different $\rho$. If we further set $\Gamma = \frac{\Vert M_R \Vert}{\Vert A  \Vert}$ and use $\Vert M_R(\xbar_k - \xbar_0) \Vert \leq \Vert M_R \Vert (\Vert \xbar_k \Vert + \Vert \xbar_0 \Vert)$, then our bound becomes 
\begin{equation*}
   \Vert b - A  \xbar_k  \Vert \lesssim \frac{1.3 c(n,k) u}{1 - \rho} \left( \Vert b \Vert + \Vert A \Vert \left( \Vert \xbar_k \Vert + \Vert \xbar_0  \Vert \right) \left( 1+ \Gamma \Vert A \Vert \Vert \Zbar_k \Vert  \right) \right).
\end{equation*}
The main aspect in which this bound differs from (3.32) in \cite{arioli2009using} is that in \cite{arioli2009using} the term $\Gamma \Vert A  \Vert \Vert \Zbar_k \Vert $ is controlled by a factor depending on $u_R$ and the precision in which the LU decomposition used as $M_R$ is computed. This comes from substitutions that rely on the specific $M_R$ when bounding $\Vert \ybar_k \Vert$. Thus, when more information on $M_R$ is available, reworking the bound for $\Vert \ybar_k \Vert$ may result in an improved bound. 

\subsection{Choosing the precisions}\label{sec:precisions}
We provide guidance on how the precisions should be set when the target backward error is of order $u$. In our experiments we observe that the achievable backward error is determined by $u + u_A \psi_A + u_L \psi_L$ and we hence ignore the term $\Vert \Zbar_k \Vert \Vert M_R(\xbar_k - \xbar_0) \Vert$ in this section. We also note that because of the structure of the former term, we do not expect the backward error to be reduced by setting $u_A$ or $u_L$ so that $u_A \psi_A \ll u$ or $u_L \psi_L \ll u$. The aim is thus to have $u \approx u_A \psi_A \approx u_L \psi_L$ in \eqref{eq:zeta2_def}. 
\begin{itemize}
    \item $u_A$. The precision for computations with $A$ should be chosen so that $u_A \approx u / \psi_A$. Numerical experiments with left-preconditioned GMRES in \cite{vieuble_thesis} show that for large $\kappa(A)$ and $\kappa(M_L)$ the quantity $\psi_A$ can be large and is driven by $\kappa(M_L)$. In such situations $u_A \ll u$ may be required. If, on the other hand, $\kappa(M_L)$ is small, then setting $u_A \approx u$ may be sufficient.  
    \item $u_L$. Guidance for setting $u_L$ comes from balancing $u \approx u_A \psi_A \approx u_L \psi_L$ and $u \approx u_L \Vert E_L M_L \Vert$. Based on the first expression, $u_L \approx u_A \psi_A / \psi_L$. Vieubl\'{e} argues that $\psi_L \leq \psi_A$ is likely, and if $\kappa(A)$ and $\kappa(M_L)$ are large then we may observe $\psi_L \ll \psi_A$ \cite{vieuble_thesis}. In these cases we can set $u_L \geq u_A$ and $u_L \gg u_A$, respectively. The quantity $\Vert E_L M_L \Vert$ depends on $M_L$ and the error in computing matrix-vector products with $M_L^{-1}$, which may be large for an ill-conditioned $M_L$. In this case thus we may require $u_L \approx u$, which is consistent with the guidance for setting $u_A \ll u$.
    \item $u_R$. Our insight on $u_R$ comes from the 
    condition \eqref{eq:rho_def_assump} (see Sections~\ref{sec:numerics_dense} - \ref{sec:numerics_sparse} for examples). It requires that $\Vert M_R \Vert \Vert E_R \Vert \ll u_R^{-1} $. Numerical experiments show that condition \eqref{eq:rho_def_assump} is sufficient but not necessary; note that a similar condition appears in \cite{arioli2009using} and it is needed to express $\Vert  \bar{y}_k \Vert $ via $\Vert 
 \xbar_k \Vert $. We can obtain a less restrictive condition on $u_R$ by keeping $\Vert \ybar_k \Vert$ but replacing $\Vert \Zbar_k \Vert$ in \eqref{eq:residual_norm_with_yk}. Using the triangle inequality in \eqref{eq:Z_computed} to bound $\Vert \Zbar_k \Vert$ we obtain terms $(u+ u_A \psi_A + u_L \psi_L) \Vert \Atilde \Vert \left( \Vert M_R^{-1} \Vert + u_R \Vert E_R \Vert \right) \Vert \ybar_k \Vert$. Thus as long as 
\[ \frac{\Vert E_R \Vert}{\Vert M_R^{-1} \Vert} \leq  u_R^{-1}, \] 
the choice of $u_R$ should not limit the backward error. $E_R$ depends on the forward error of matrix-vector products with $M_R$. If $\kappa(M_R)$ is large, we may need a small $u_R$ for the condition to be satisfied. When $\Vert M_R \Vert $ and $\kappa(M_R)$ are small, a large value for $u_R$ may suffice. Note that these comments take into account the backward error only and not the FGMRES iteration count.
\end{itemize}
  
\subsection{Forward error}\label{sec:forward_error}
A rule of thumb says that the forward error can be bounded by multiplying the backward error by the condition number of the coefficient matrix; see, for example, \cite{high:ASNA2}. 
Using \eqref{eq:backward_error} thus gives the bound
\begin{equation}\label{eq:forw_err_left_precond}
    \frac{\Vert x - \xbar_k \Vert}{\Vert x \Vert} \leq  \frac{1.3 c(n,k)}{1 - \rho}\zeta \kappa(\Atilde),
\end{equation}
where $x$ is the solution to \eqref{eq:Ax=b_split_precond} and $\xbar_k$ is the output of FGMRES. 
Note that the bound depends on the condition number of the left-preconditioned matrix $\Atilde$. If $\kappa(\Atilde)$ and $\zeta$ are small, then the forward error is small too, and thus $\xbar_k \approx x$. Then $b - A\xbar_k = A(x - \xbar_k)$ is small and implies a small backward error with respect to the original system as previously noted.

The forward error bound can also be formulated with respect to the split-preconditioned matrix $\Ahat \coloneqq M_L^{-1} A M_R^{-1}$ as follows
\begin{equation}\label{eq:forw_err_split_prec}
    \frac{\Vert x - \xbar_k \Vert}{\Vert x \Vert} \leq  \frac{1.3c(n,k)}{1 - \rho} \zeta \kappa(\Ahat) \kappa(M_R).
\end{equation}
Note that \eqref{eq:forw_err_split_prec} is weaker than \eqref{eq:forw_err_left_precond} as $\kappa(\Atilde) \leq \kappa(\Ahat) \kappa(M_R)$. However, \eqref{eq:forw_err_split_prec} may be useful if $\kappa(\Ahat)$ and $ \kappa(M_R)$ or their estimates are known and such information is not available for $\kappa(\Atilde)$. 
The bounds \eqref{eq:forw_err_left_precond} and \eqref{eq:forw_err_split_prec} suggest that guaranteeing a small forward error requires controlling the backward error and constructing the preconditioners so that either $\kappa(\Atilde)$ or both $\kappa(\Ahat)$ and $ \kappa(M_R)$ are small (depending on which condition numbers can be evaluated). If $A$ is ill-conditioned, then achieving a small $\kappa(\Atilde)$ in \eqref{eq:forw_err_left_precond} requires an $M_L$ with a high condition number. Note that in this case, as discussed in the previous section, we may have to set $u_A \ll u$ and can get away with $u_L \gg u_A$. The bound \eqref{eq:forw_err_split_prec} indicates that if we achieve a small $\kappa(\Ahat)$ at the price of $\kappa(M_R) \approx \kappa(A)$ then we cannot guarantee a smaller forward error than when no preconditioning is used because unpreconditioned FGMRES is equivalent to unpreconditioned GMRES in uniform precision with backward error bounded by $\frac{cnu}{1-cnu}$, where $c$ is a constant \cite{paige2006modified}.

\subsection{Left-, right-, or split-preconditioning}
As previously stated, the split-preconditioning approach allows us to analyze the left- and right-preconditioned cases as well. In this section, we explore which preconditioning strategy may be preferred under certain objectives. The discussion is based on the bounds for the backward and forward error. We first simplify these for left- and right-preconditioning.

If only left-preconditioning is used, $M_R = I$, $E_R = 0$, $\Zbar_k = \Vbar_k$, and $\rho = 1.3c(n,k)u$, and Algorithm~\ref{alg:fgmres} is equivalent to left-preconditioned GMRES. The dominant term in the relative backward error for the preconditioned system is thus 
\begin{equation}\label{eq:zeta2_left-only_precond}
    \zeta_2 = \left(u + u_A \psi_A + u_L \psi_L \right) \Vert \Atilde \Vert \left( \Vert \xbar_k - \xbar_0 \Vert + \Vert \xbar_0  \Vert \right) \approx  \left(u + u_A \psi_A + u_L \psi_L \right) \Vert \Atilde \Vert \Vert \xbar_k \Vert,
\end{equation}
where the approximation holds if $\Vert \xbar_k \Vert \approx  \Vert \xbar_0 \Vert$ or $\xbar_0 =0$. This is equivalent to the result for left-preconditioned GMRES in \cite[Theorem 7.1]{vieuble_thesis}. 

In the right-preconditioning case, $M_L = I$, $E_L = 0$ and $\psi_L = 0$, and $u_A \psi_A = \max_j \{ \Vert \Delta A_j z_j\Vert / \Vert  A \Vert \Vert z_j\Vert \}$. Then \eqref{eq:backward_error} gives the bound for the relative backward error for the original system \eqref{eq:Ax=b} and
\begin{gather*}
    \zeta_1 = u \Vert b \Vert, \text{ and} \\
    \zeta_2 = \left( u  + u_A \psi_A  \right) \Vert A \Vert  \left( \Vert \Zbar_k \Vert \Vert M_R(\xbar_k - \xbar_0) \Vert + \Vert \xbar_0 \Vert \right).
\end{gather*}
We make the following observations.
\begin{itemize}
    \item Consider the case where a small backward error is the main concern and $A$ is ill-conditioned. If we have a `good' preconditioner, so that $\kappa(\Atilde)$ is small and we can afford setting $u_A$ and $u_L$ to precisions that are high enough to neutralize the $\psi_A$ and $\psi_L$ terms, then Corollary~\ref{cor:be_original_system} and \eqref{eq:zeta2_left-only_precond} can guarantee a small relative backward error when full left-preconditioning is used. If however, we cannot afford setting $u_A$ and $u_L$ to high precisions but can construct a split-preconditioner such that $\kappa(M_L)$ is small, then split-preconditioning (note that in this case $\psi_A$ and $\psi_L$ may be smaller too) or full right-preconditioning may be preferential. Note however that a small backward error for these options can only be guaranteed if $\Vert \Zbar_k \Vert$ is expected to be small or if the bounds for $\Vert \ybar_k \Vert$ are reworked taking into account a specific preconditioner in order to control $\Vert \Zbar_k \Vert$ (as for the right-preconditioning with an LU decomposition in \cite{arioli2009using}).
    \item If we are aiming for a small forward error and can obtain a small relative backward error with full left-preconditioning as detailed above, then this approach also gives a small forward error. If we, however, do not have the flexibility of setting $u_A$ and $u_L$, then it is not clear from the bounds which preconditioning approach gives the best results.
    \item Assume that our main concern is applying the preconditioner in lower than the working precision. This may be relevant, for example, when $A$ is very sparse and the preconditioner uses some dense factors. In this case, the bounds suggest that full left-preconditioning should not be used as $u_A \psi_A$ and $u_L \psi_L$ may be large. Full right-preconditioning may be suitable in this case although the bound is affected by $\Vert \Zbar_k \Vert$.
\end{itemize}

\section{Example: LU preconditioner}\label{sec:example_LU_precond}
We supplement the theoretical analysis in the previous section with an example. Assume that an approximate LU decomposition of $A$ is computed, for example in a low precision, and the computed factors $\Lbar$ and $\Ubar$ are used for preconditioning. We choose this preconditioning due to its effectiveness and ease of application; note that there is no structural advantage to applying it as a left-, split- or right-preconditioner. 

If split-preconditioning is used, then $M_L = \Lbar$ and $M_R = \Ubar$. 
In Algorithm~\ref{alg:fgmres}, products with $A$ are computed in precision $u_A$ and hence
\begin{equation}\label{eq:bound_psiA}
    \psi_{A,j} = \frac{\Vert M_{L}^{-1} \Delta A_j z_j\Vert}{u_A \Vert \Atilde \Vert \Vert z_j \Vert} \leq \ctilde_1(n) \frac{ \Vert \lvert \Lbar^{-1} \rvert \lvert A \rvert \Vert \Vert z_j \Vert }{ \Vert \Lbar^{-1} A \Vert  \Vert z_j \Vert} = \ctilde_1(n) \frac{ \Vert \lvert \Lbar^{-1} \rvert \lvert A \rvert \Vert }{ \Vert \Lbar^{-1} A \Vert },
\end{equation}
where $\ctilde_i(n)$ is a constant that depend on $n$. We expect $\Vert \lvert \Lbar^{-1} \rvert \lvert A \rvert \Vert / \Vert \Lbar^{-1} A \Vert$ to be moderate for many systems and in this case setting $u_A = u$ may be sufficient.

We apply $M_L$ by solving a triangular system $\Lbar w_j = (A +  \Delta A_j) z_j$ via substitution in precision $u_L$. From standard results we know that the computed $\wbar_j$ satisfies
\begin{gather*}
    (\Lbar + \Delta L_j) \wbar_j = (A +  \Delta A_j) z_j, \\
    \text{where } \lvert  \Delta L_j \rvert \leq \ctilde_2(n) u_L \lvert  \Lbar \rvert.
\end{gather*}
Thus 
\begin{equation*}
\Delta M_{L,j} = \Lbar^{-1} - (\Lbar + \Delta L_j)^{-1} \approx \Lbar^{-1} \Delta L_j \Lbar^{-1}.    
\end{equation*}
We use this to bound $\psi_{L,j}$ as
\begin{equation}\label{eq:bound_psiL}
    \psi_{L,j} =  \frac{\Vert \Delta M_{L,j} A z_j\Vert}{u_L \Vert \Atilde \Vert \Vert z_j \Vert} \approx \frac{\Vert \Lbar^{-1} \Delta L_j \Lbar^{-1} A z_j\Vert}{u_L \Vert \Lbar^{-1} A \Vert \Vert z_j \Vert} \leq  \frac{\Vert \Lbar^{-1}   \Delta L_j \Vert \Vert \Lbar^{-1}  A \Vert \Vert z_j \Vert}{u_L \Vert \Lbar^{-1} A \Vert \Vert z_j \Vert}  \leq  \ctilde_3(n) \kappa_2(\Lbar).
\end{equation}
The bound \eqref{eq:bound_psiL} is obtained using the bound on the forward error of solving a triangular system. In general such systems are solved to high accuracy and thus we expect \eqref{eq:bound_psiL} to be a large overestimate. Note that bounds \eqref{eq:bound_psiA} and \eqref{eq:bound_psiL} hold for every $j$.

When using full left-preconditioning, $M_L = \Lbar \Ubar $ and $M_R = I $. This case is considered in \cite{vieuble_thesis}. We bound $\psi_{A,j}$ in a similar way as in the split-preconditioned case, i.e., 
\begin{equation*}
    \psi_{A,j} \leq \ctilde_4(n) \frac{ \Vert \lvert \Ubar^{-1}  \Lbar^{-1} \rvert \lvert A \rvert \Vert  }{ \Vert \Ubar^{-1} \Lbar^{-1} A \Vert},
\end{equation*}
which can be expected to be moderate in many cases as well.
Applying $M_L$ now requires consecutively solving two triangular systems of equations, that is, 
\begin{gather*}
    (\Lbar + \Delta L_j) (\Ubar + \Delta U_j) \wbar_j = (A +  \Delta A_j) z_j, \\
    \text{where}\quad \lvert  \Delta U_j \rvert \leq \ctilde_5(n) u_L \lvert  \Ubar \rvert \quad\text{and}\quad \lvert  \Delta L_j \rvert \leq \ctilde_6(n) u_L \lvert  \Lbar \rvert.
\end{gather*}
This gives
\begin{equation}\label{eq:deltaM_L_full_left-precond}
\Delta M_{L,j} = \Ubar^{-1}\Lbar^{-1} - (\Ubar + \Delta U_j)^{-1} (\Lbar + \Delta L_j)^{-1} \approx  \Ubar^{-1} \Lbar^{-1} \Delta L_j \Lbar^{-1} + \Ubar^{-1} \Delta U_j \Ubar^{-1} \Lbar^{-1},    
\end{equation}
where we omit the term involving $\Delta L_j \Delta U_j $, and hence
\begin{align*}
    \psi_{L,j}  \leq & \frac{\Vert  \Ubar^{-1} \Lbar^{-1} \Delta L_j \Lbar^{-1}  A z_j\Vert + \Vert \Ubar^{-1} \Delta U_j \Ubar^{-1} \Lbar^{-1}  A z_j\Vert }{u_L \Vert \Ubar^{-1} \Lbar^{-1} A \Vert \Vert z_j \Vert} \\
    \leq & \frac{\Vert  \Ubar^{-1} \Vert \Vert \vert \Lbar^{-1} \vert \vert L \vert \Vert \Vert \Lbar^{-1}  A \Vert }{\Vert \Ubar^{-1} \Lbar^{-1} A \Vert } + \Vert \vert \Ubar^{-1}\vert \vert  U \vert  \Vert \\
    = & \frac{\Vert  \Ubar^{-1} \Vert  \Vert \Lbar^{-1}  A \Vert }{\Vert \Ubar^{-1} \Lbar^{-1} A \Vert }\textrm{cond}(\Lbar) + \textrm{cond}(\Ubar),
\end{align*}
where $\textrm{cond}(B) = \Vert \vert B^{-1} \vert \vert B \vert \Vert$. As discussed in \cite{vieuble_thesis}, $\textrm{cond}(\Lbar) $ and $ \textrm{cond}(\Ubar)$ are modest when using partial pivoting. The term $\Vert  \Ubar^{-1} \Vert $, however, can be large for an ill-conditioned $A$.

For full right-preconditioning, $M_L = I $ and $M_R = \Lbar \Ubar $. Then $ \psi_{L,j} = 0$ and
\begin{equation}\label{eq:bound_psiA_right_precond}
    \psi_{A,j} = \frac{\Vert \Delta A_j z_j\Vert}{u_A \Vert A \Vert \Vert z_j \Vert} \leq \frac{\Vert \vert A \vert \Vert}{\Vert A \Vert } \leq \sqrt{n},
\end{equation}
where the final inequality is due to $\Vert \vert A \vert \Vert_2 \leq \sqrt{\textrm{rank}(A)} \Vert A \Vert_2 $ \cite[Lemma 6.6]{high:ASNA2}. We can thus use $u_A = u$ in most cases.

\subsection{Numerical example: synthetic dense systems with split-preconditioning}\label{sec:numerics_dense}
We perform numerical experiments in MATLAB R2021a\footnote{The code is available at \url{https://github.com/dauzickaite/mpfgmres/}.} using a setup similar to an example from \cite{arioli2009using}. An $n \times n$ coefficient matrix $A=UDV$ is constructed by generating random orthogonal $n \times n$ matrices $U$ and $V,$ and setting $D$ to be diagonal with elements $10^{-c(j-1)/(n-1)}$ for $j=1,2,\dots,n$. The condition number of $A$ is $10^c$ and we vary its value. The right hand side $b$ is a random vector with uniformly distributed entries. The preconditioner is computed as a low precision LU factorization. Namely, for $c \in \{1,2,\dots,5 \}$ we use $[\Lbar,\Ubar] = lu(mp(A,4))$, where $mp(\cdot,4)$ calls the Advanpix Multiprecision Computing Toolbox \cite{advanpix} and simulates precision accurate to four decimal digits; note that this has a smaller unit roundoff than IEEE half precision (see Table~\ref{tab:ieee_param} for the unit roundoff values). For $c \in \{6,7,\dots,10 \}$ we compute LU factorization in single precision using the built-in MATLAB single precision data type. We set $M_L = \Lbar$ and $M_R = \Ubar$, and $E_L = \vert \Lbar^{-1}\vert  \vert \Lbar \vert \vert \Lbar^{-1}\vert $ and $E_R = \vert \Ubar^{-1} \vert  \vert \Ubar \vert \vert \Ubar^{-1}\vert$. The left-preconditioner can slightly reduce the condition number of the coefficient matrix whereas the split-preconditioner achieves a high reduction (Table~\ref{tab:cond_no_different_c}). 

\begin{table}[]
 \caption{Unit roundoff $u$ for IEEE floating point arithmetics.}
    \centering
    \begin{tabular}{l|c}
        Arithmetic &  $u$  \\
         \hline
        fp16 (half)  &  $2^{-11} \approx 4.88 \times 10^{-4}$ \\
        fp32 (single) & $2^{-24} \approx 5.96 \times 10^{-8}$ \\
        fp64 (double) & $2^{-53} \approx 1.11 \times 10^{-16}$ \\
        fp128 (quadruple) & $2^{-113} \approx 9.63 \times 10^{-35}$
    \end{tabular}
    \label{tab:ieee_param}
\end{table}

We set the working precision $u$ to double. Bounds for $\psi_A$ in Table~\ref{tab:cond_no_different_c} indicate that there is no need for $u_A < u$, thus we choose $u_A = u$. The preconditioners are applied using all combinations of half, single, double, and quadruple precisions. Half precision is simulated via the \textit{chop} function \cite{HighamChop}, and Advanpix is used for quadruple precision. We expect $\kappa(M_L)$ to be a large overestimate for $\psi_L$. $\kappa(M_R)$ suggests that the condition $\rho <1$ in \eqref{eq:rho_def_assump} should be satisfied with $u_R$ set to any of the four precisions, except half for large $c$ values. We approximate $\Vert E_R \Vert / \Vert M_R^{-1} \Vert$ by $\min\{\Vert \vert \Ubar^{-1} \vert  \vert \Ubar \vert \Vert, \Vert \vert \Ubar \vert \vert \Ubar^{-1}\vert \Vert \}$ in Table~\ref{tab:cond_no_different_c}, where we round to the nearest whole number. 
If $u_R$ is set to half precision, then $u_R^{-1} = 2048$, which is slightly smaller than the values of $\min\{\Vert \vert \Ubar^{-1} \vert  \vert \Ubar \vert \Vert, \Vert \vert \Ubar \vert \vert \Ubar^{-1}\vert \Vert \}$ for $c \geq 6$. This indicates that applying $M_R$ in half precision may affect the backward error, however note that our choice for $E_R$ is expected to be an overestimate. The solver tolerance $\tau$ (see Algorithm~\ref{alg:fgmres}) is set to $4u$ and we use $x_0 = 0$. For the unpreconditioned system, FGMRES converges in 200 iterations when $c=1$ and does not converge for other $c$ values.

\begin{figure}
    \centering
\begin{subfigure}[t]{0.45\linewidth}
  \centering
 \includegraphics[width=\linewidth]{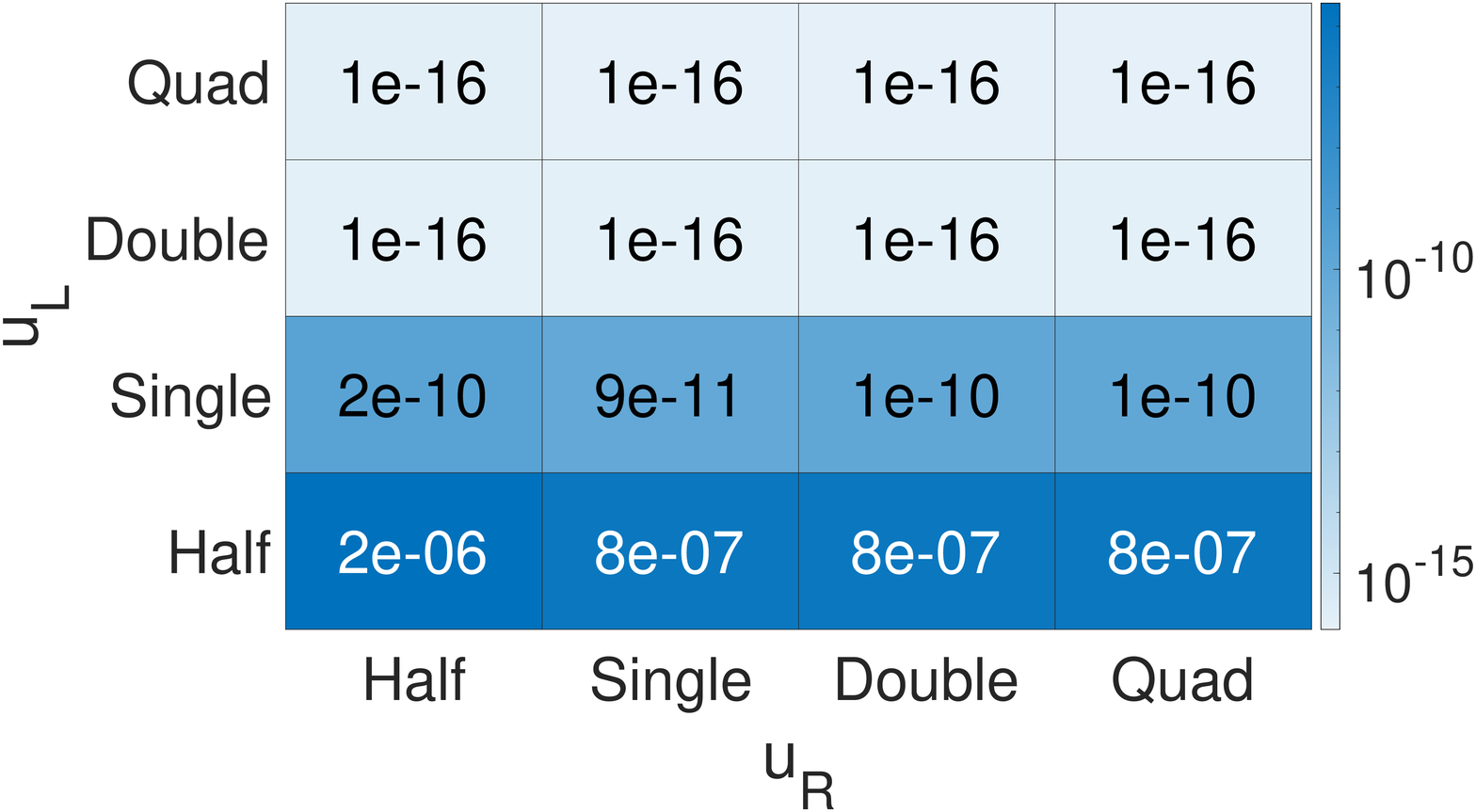}
  \caption{BE}
\end{subfigure}
\begin{subfigure}[t]{0.45\linewidth}
  \centering
 \includegraphics[width=\linewidth]{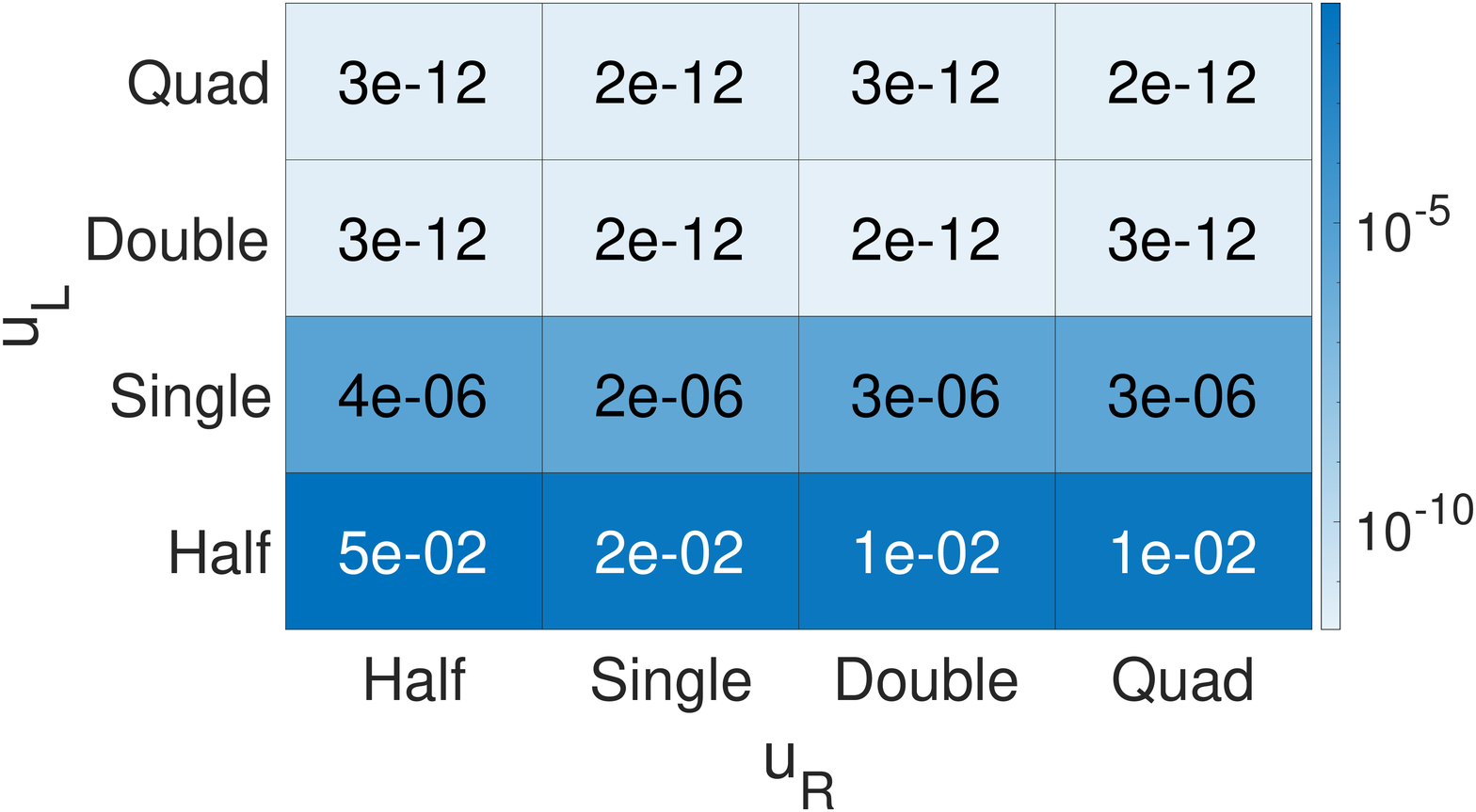}
  \caption{FE}
\end{subfigure}
\begin{subfigure}[t]{0.45\linewidth}
  \centering
 \includegraphics[width=\linewidth]{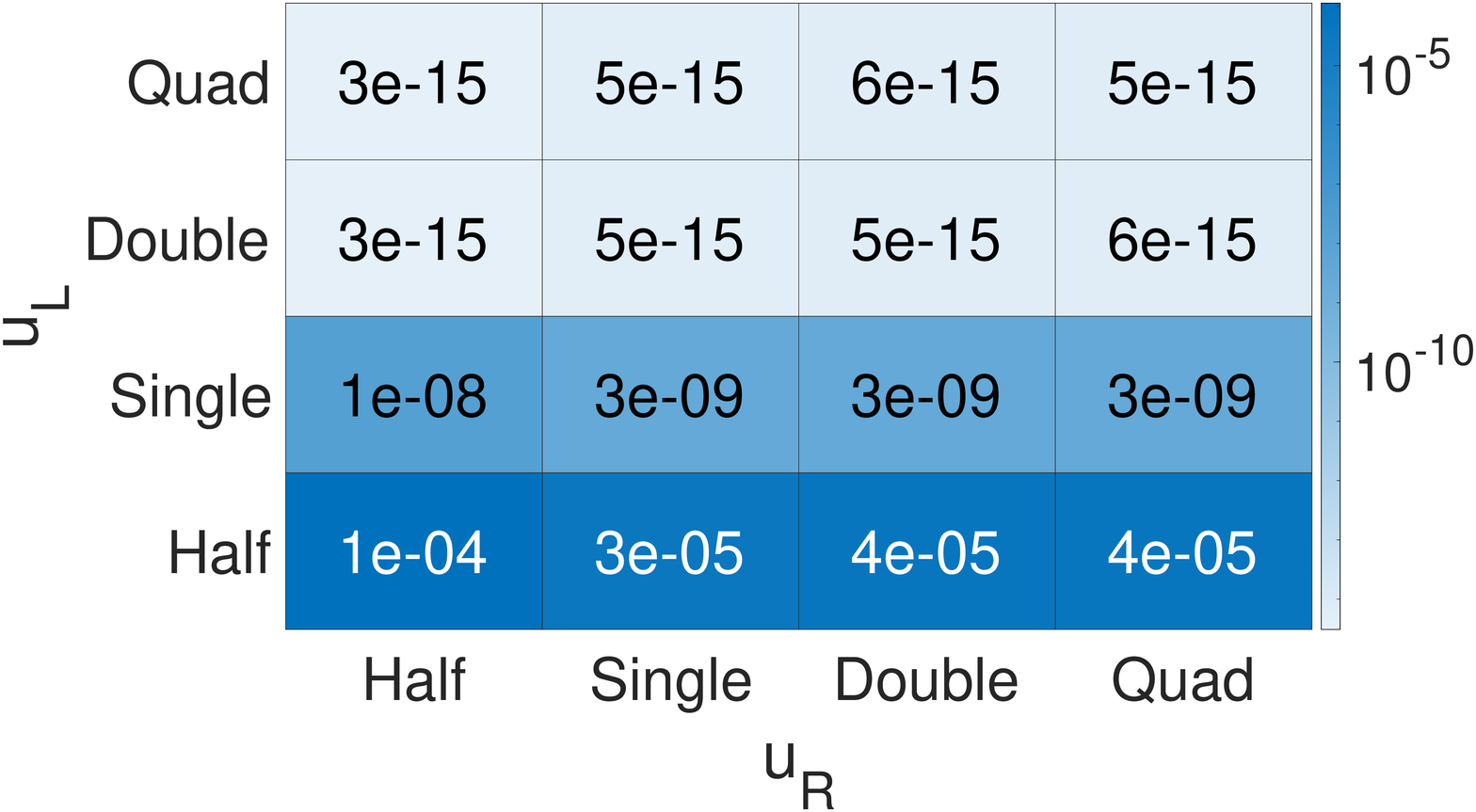}
  \caption{$\zeta$}
\end{subfigure}
\begin{subfigure}[t]{0.45\linewidth}
  \centering
 \includegraphics[width=\linewidth]{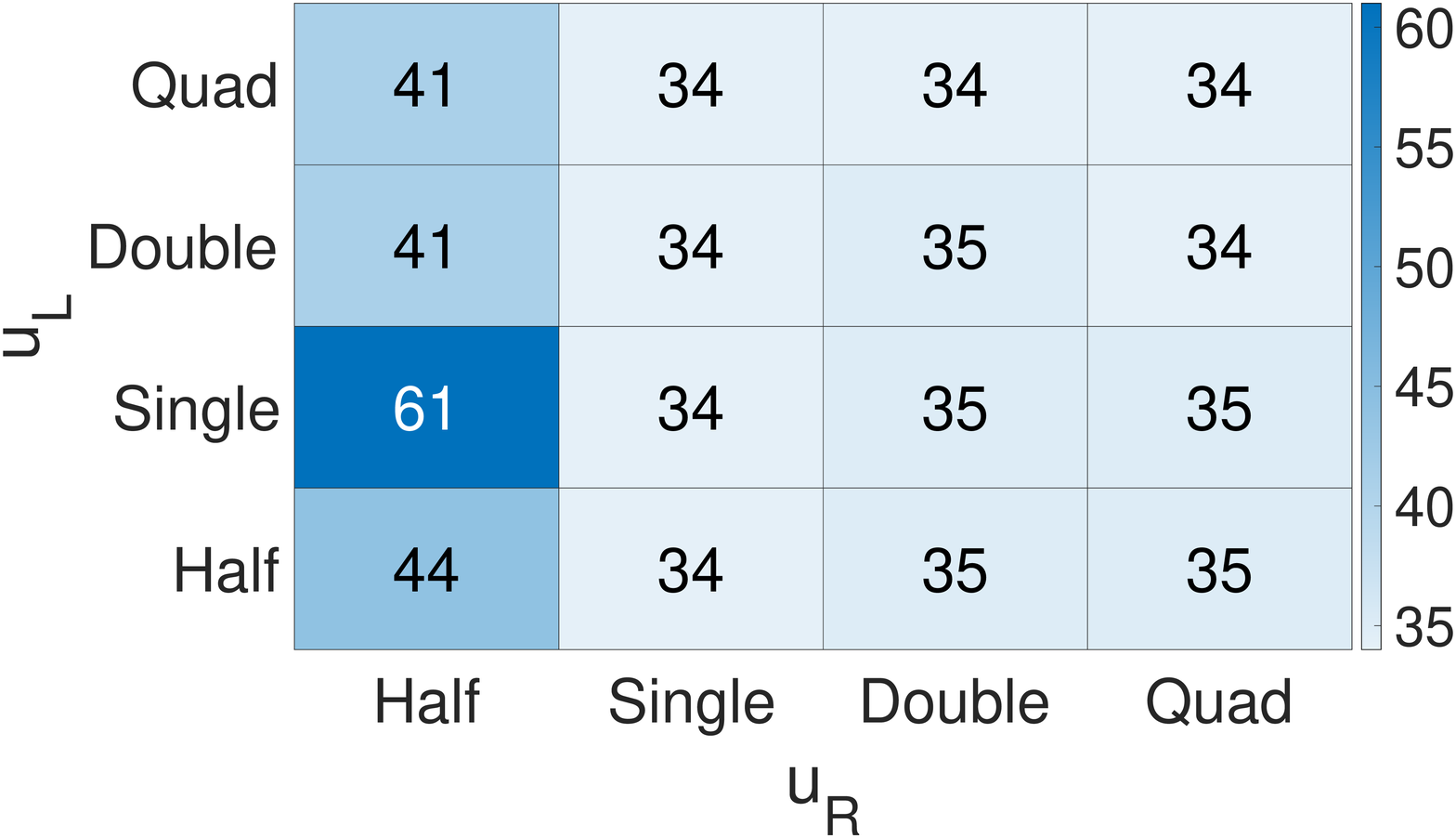}
  \caption{iterations}
\end{subfigure}
\begin{subfigure}[t]{0.45\linewidth}
  \centering
 \includegraphics[trim={0 0 0 3},clip,width=\linewidth]{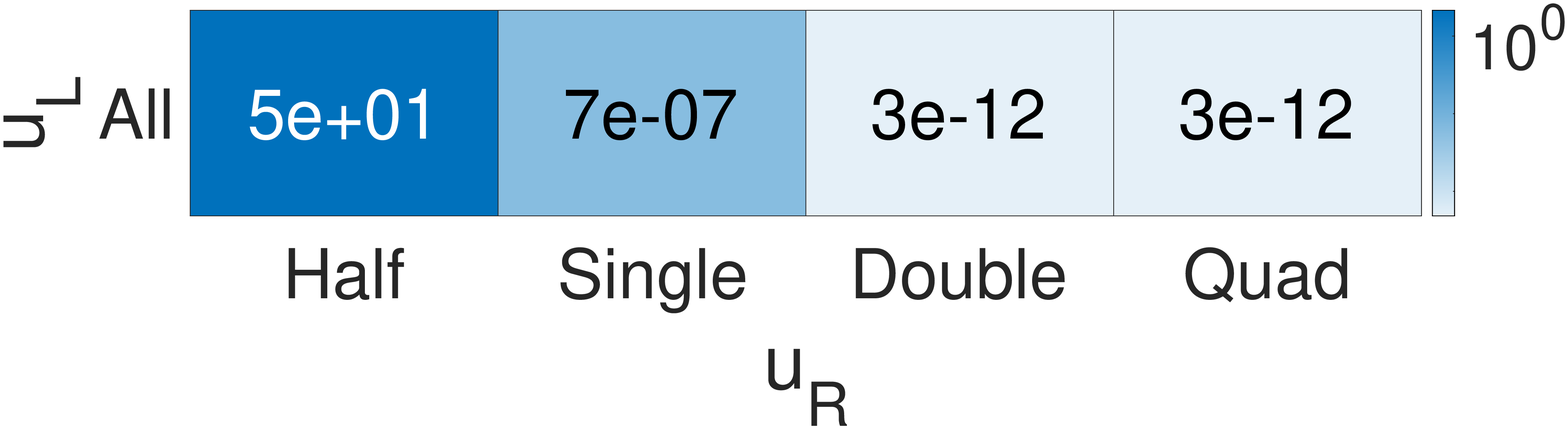}
  \caption{$\rho$}
\end{subfigure}
    \caption{Synthetic problem, $c=5$. BE is the relative backward error and FE is the relative forward error, $\zeta$ is as defined in \eqref{eq:zeta}, and $\rho$ is as defined in \eqref{eq:rho_def_assump}.}
    \label{fig:dense_c5}
\end{figure}

We show results for $c=5$ for all precision combinations (Figure~\ref{fig:dense_c5}), and for all $c$ values with $u_L$ set to single and $u_R$ set to double (Table~\ref{tab:results_uL_single}), and $u_L$ set to double and $u_R$ set to single (Table~\ref{tab:results_uR_single}). We report the relative backward error (BE) of the original problem, that is,
\begin{equation*}
    \frac{\Vert b - A \xbar_k \Vert}{\Vert b \Vert + \Vert A \Vert \Vert \xbar_k \Vert}
\end{equation*}
and compute the dominant part of the backward error bound $\zeta$ (as defined in \eqref{eq:zeta}). Note that $\zeta$ bounds the relative backward error for the left-preconditioned system.

From Figure~\ref{fig:dense_c5}, we can see that the achievable backward error and subsequently the forward error depends on $u_L$. As expected from theory, $u_R$ does not affect the achievable backward error, however $u_R$ influences the iteration count. Setting $u_L$ to half results in extra iterations when $c=1$, $c=2$ and $c=6$ (not shown). Note that setting $u_L$ to quadruple, and $u_R$ to double or quadruple does not give any benefit. As mentioned, the backward error bound \eqref{eq:backward_error} is dominated by $\zeta_2$. From Tables~\ref{tab:results_uL_single} and \ref{tab:results_uR_single}, we can see that the quantity $\Vert \Zbar_k \Vert \Vert M_R(\xbar_k - \xbar_0) \Vert$ can become large, but it stays of the order of $\Vert \xbar_k \Vert$ or close to it (not shown) and thus $\zeta$ gives a good estimate of the backward error. If, however, $\Vert x \Vert$ is small, then the term $\Vert \Zbar_k \Vert \Vert M_R(\xbar_k - \xbar_0) \Vert$ can impair the bound. The increase in the forward error compared to the backward error is well estimated by $\kappa(\Atilde)$, whereas $\kappa(\Ahat) \kappa(M_R)$ is an overestimate. From Figure~\ref{fig:dense_c5}, we see that the condition $\rho <1$ is sufficient, but not necessary as it is not satisfied when $c=5$ and $u_R$ is set to half.

\begin{table}[]
    \caption{Synthetic problems. Condition numbers of the unpreconditioned and preconditioned coefficient matrices and preconditioners, and bounds for $\psi_A$ and $\psi_L$ ($\kappa(M_L)$ is also the bound for $\psi_L$). The preconditioners are computed in precision accurate to four decimal digits for $c<6$ and in single precision for $c\geq 6$.}
    \centering
    \begin{tabular}{c|c|c|c|c|c|c|c}
c & $\kappa(A)$	& $\kappa(\Atilde)$	&$\kappa(\Ahat)$ & $\kappa(M_R)$ & $\kappa(M_L)$	& $\psi_A$ bound &  $\frac{\Vert E_R \Vert }{ \Vert M_R^{-1} \Vert} $  \\	
& & & & & & & approx. \\
\hline
1 & 10&	$2.69 \times 10^2$&	1.06 &	$2.69 \times 10^2$	& $3.05 \times 10^2$ &	$4.85 \times 10^1$ & 615\\
 2 & $10^2$	& $4.80\times 10^2$	& 1.15	& $4.80 \times 10^2$	& $3.23 \times 10^2$	& $6.30 \times 10^1$ & 1551 \\
 3 & $10^3$& 	$1.69 \times 10^3$& 	1.66	& $1.69 \times 10^3$ &	$3.09 \times 10^2$	& $7.27 \times 10^1$ & 1853 \\
 4 & $10^4$	& $1.40 \times 10^4$ & 	$1.06 \times 10^1$	& $1.43 \times 10^4$	& $2.79 \times 10^2$	& $1.64 \times 10^2$ & 1949\\
 5 & $10^5$	& $7.73 \times 10^4$ &	$3.25 \times 10^2$	& $1.10 \times 10^5$	& $2.85 \times 10^2$	& $1.18 \times 10^2$ & 2048 \\
 \hline
 6 & $10^6$ &	$4.68 \times 10^5$	& $1.14$	& $4.69 \times 10^5$	& $2.64 \times 10^2$& 	$2.54 \times 10^2$ & 2144\\
 7 & $10^7$ & $2.76 \times 10^6$ & $2.34$ & $2.74 \times 10^6$ & $3.50 \times 10^2$ & $1.80 \times 10^2$ & 2336\\
 8 & $10^8$ & $3.98 \times 10^7$ & $7.00 \times 10^1$ & $4.04 \times 10^7$ & $3.41 \times 10^2$ & $1.03 \times 10^2$ & 2551\\
9 & $10^9$ & $3.19  \times 10^8$ & $5.44 \times 10^3$ & $4.78 \times 10^8$ & $3.95 \times 10^2$ & $1.37 \times 10^2$ & 2716 \\
10 & $10^{10}$ & $4.49 \times 10^9$ & $5.88 \times 10^4$ & $5.72 \times 10^8$ & $3.06 \times 10^2$ & $1.46 \times 10^2$ & 2768
    \end{tabular}
    \label{tab:cond_no_different_c}
\end{table}

\begin{table}[]
    \caption{Synthetic problems with $u_L$ single and $u_R$ double. IC denotes the iteration count, BE is the relative backward error and FE is the relative forward error. For $c=10$, the solver is terminated at 200 iterations without satisfying the convergence criteria.}
    \centering
     \scalebox{0.7}{
    \begin{tabular}{c|c|c|c|c|c|c|c|c}
c	& IC & BE &	FE	&	$\zeta$&	$\Vert \Zbar_k \Vert \Vert M_R(\xbar_k - \xbar_0) \Vert$	& $\psi_A$&	$\psi_L$	& $\rho$ \\
\hline
1 & 6 & $2.45 \times 10^{-7}$& 	$1.61 \times 10^{-6}$	&$2.91 \times 10^{-7}$& $ 5.04 \times 10^1$	& 1.21	& 2.22 &	$1.44 \times 10^{-15}$ \\
2	& 7 & $3.82 \times 10^{-8}$	& $1.31 \times 10^{-6}$ & $9.34 \times 10^{-8}$	& $2.92 \times 10^2$ 	& 2.08	& $7.46 \times 10^{-1}$	& $3.59 \times 10^{-15}$\\
3	& 9 & $5.76 \times 10^{-9}$	& $1.46 \times 10^{-6}$		& $ 2.56 \times 10^{-8}$&	$2.53 \times 10^3$&	3.79	&$1.88  \times 10^{-1}$	& $1.51 \times 10^{-14}$ \\
4	& 15 &$5.45 \times 10^{-10}$	& $1.27 \times 10^{-6}$	&	$ 3.14 \times 10^{-9}$	& $ 2.30 \times 10^4$	& 3.55	& $2.19  \times 10^{-2}$	& $1.46 \times 10^{-13}$\\
5 & 	35 & $1.03 \times 10^{-10}$&	$2.59 \times 10^{-6}$	& $ 2.51 \times 10^{-9}$&	$  8.39 \times 10^5$	&4.58	&$4.28  \times 10^{-3}$ &	$2.66 \times 10^{-12}$\\
 \hline
6	& 7 & $6.37 \times 10^{-12}$	& $1.58 \times 10^{-6}$	&  $ 5.30 \times 10^{-11}$	& $ 1.35 \times 10^6$& 	4.81	& $5.08 \times 10^{-4}$	& $4.88 \times 10^{-12}$ \\
7 & 11 & $5.88  \times 10^{-13}$ & $9.77  \times 10^{-7}$ & $ 8.25 \times 10^{-12}$ & $1.42 \times 10^7$ & 8.09 & $6.88  \times 10^{-5}$ & $3.28 \times 10^{-11}$\\
8 & 21 & $7.11 \times 10^{-14}$ & $8.24 \times 10^{-7}$ & $1.72 \times 10^{-12}$ & $2.71 \times 10^8$ & 6.08 & $6.95 \times 10^{-6}$ & $8.31 \times 10^{-10}$ \\
9 & 92 & $ 7.97 \times 10^{-14}$ & $8.82 \times 10^{-6}$ & $1.32 \times 10^{-11}$ & $4.02 \times 10^{10}$ & 6.82 & $3.33 \times  10^{-6}$ & $4.48  \times 10^{-8}$ \\
10 & 200 & $1.21 \times 10^{-13} $ & $2.09 \times 10^{-4}$ & $2.30 \times 10^{-11}$ & $5.95 \times 10^{11}$ & 8.28 & $3.68 \times  10^{-6}$ & $6.97 \times 10^{-8}$
    \end{tabular}
    }
    \label{tab:results_uL_single}
\end{table}

\begin{table}[]
    \caption{Synthetic problems with $u_L$ double and $u_R$ single. IC denotes the iteration count, BE is the relative backward error and FE is the relative forward error. For $c=10$, the solver is terminated at 200 iterations without satisfying the convergence criteria.}
    \centering
     \scalebox{0.7}{
    \begin{tabular}{c|c|c|c|c|c|c|c|c}
 c & IC & BE &	FE	& 	$\zeta$&	$\Vert \Zbar_k \Vert \Vert M_R(\xbar_k - \xbar_0) \Vert$	& $\psi_A$&	$\psi_L$	& $\rho$ \\
\hline
1	& 6 & $5.14 \times 10^{-16}$&	$ 3.04 \times 10^{-15}$		& $1.17\times 10^{-15}$	&$ 5.05 \times 10^1$ &	1.26	& 2.48	&$1.29 \times 10^{-9}$\\
2	& 7 &$1.20 \times 10^{-16}$	& $3.61 \times 10^{-15}$&	$ 7.69 \times 10^{-16}$&	$ 2.92 \times 10^2$ &	1.66 &	$6.33 \times 10^{-1}$&	$3.06 \times 10^{-9}$\\
3	& 9 & $8.62 \times 10^{-17}$	& $2.07 \times 10^{-14}$ & $ 1.19 \times 10^{-15}$&	$2.53 \times 10^3$	& 3.46	& $2.45 \times 10^{-1}$&	$1.31 \times 10^{-8}$ \\
4 &	15 &$6.91 \times 10^{-17}$&	$1.36 \times 10^{-13}$		& $ 1.12 \times 10^{-15}$	& $ 2.30 \times 10^4$&	3.19	& $2.44 \times 10^{-2}$	& $1.05 \times 10^{-7}$ \\
5	& 34 &  $1.42 \times 10^{-16}$	& $2.45 \times 10^{-12}$	& 	$5.20 \times 10^{-15}$	& $ 8.20\times 10^5$&	3.86	& $5.87 \times 10^{-3}$	& $6.88 \times 10^{-7}$ \\
 \hline
6	& 7 & $5.17 \times 10^{-17}$	& $8.85 \times 10^{-12}$&		$ 9.97 \times 10^{-16}$	& $1.31 \times 10^6$	& 4.28 &	$7.65 \times 10^{-4}$	&$3.81 \times 10^{-6}$ \\
7 & 11 & $4.97 \times 10^{-17}$ & $6.91 \times 10^{-11}$ & $1.59 \times 10^{-15}$ & $1.40 \times 10^7$ & 6.23 & $1.11 \times 10^{-4}$ & $2.40 \times 10^{-5}$\\
8 & 21 & $6.58 \times 10^{-17}$ & $5.29 \times 10^{-10}$ & $2.41  \times 10^{-15}$ & $2.29 \times 10^8$ & 5.18 & $1.71 \times 10^{-5}$ & $3.38 \times 10^{-4}$\\
9 & 158 & $2.91 \times 10^{-16}$ & $4.98 \times 10^{-8}$ & $7.65 \times 10^{-14}$ & $4.28 \times 10^{10}$ & 8.79 & $7.92\times 10^{-6} $ & $3.16 \times 10^{-3}$\\
10 & 200 & $4.69 \times 10^{-16}$ & $5.71 \times 10^{-7}$ & $8.81 \times 10^{-14}$ & $5.24 \times 10^{11}$ & 7.62 & $7.24 \times 10^{-6}$ & $4.38 \times 10^{-3}$
    \end{tabular}
    }
    \label{tab:results_uR_single}
\end{table}

\subsection{Numerical example: left- and right-preconditioning}
We also perform experiments with full left- (i.e., $M_L = \Lbar \Ubar$ and $M_R = I$) and full right-preconditioning (i.e., $M_L = I$ and $M_R = \Lbar \Ubar$) using the same set-up as in the previous section. For the left-preconditioning case using \eqref{eq:deltaM_L_full_left-precond} we set $E_L = \vert \Ubar^{-1} \Lbar^{-1} \vert \vert  \Lbar \vert \vert  \Lbar^{-1} \vert + \vert \Ubar^{-1}  \vert \vert \Ubar \vert \vert \Ubar^{-1} \Lbar^{-1} \vert$ and $E_R = 0$. Equivalently, for the right-preconditioning case $E_R = \vert \Ubar^{-1} \Lbar^{-1} \vert \vert  \Lbar \vert \vert  \Lbar^{-1} \vert + \vert \Ubar^{-1}  \vert \vert \Ubar \vert \vert \Ubar^{-1} \Lbar^{-1} \vert$ and $E_L = 0$. Both left- and right-preconditioners are effective in reducing the condition number of the coefficient matrix (see Table~\ref{tab:left-right-precond_condition_numbers}). Note that the bounds for $\psi_A$ and $\psi_L$ in the left-preconditioning case indicate that $\psi_L$ may have to be applied in precision such that $u_L \leq u$ for all problems and we may need $u_A \leq u$ for highly ill-conditioned cases. For right-preconditioning, from \eqref{eq:bound_psiA_right_precond} we know that $\psi_A$ can be bounded by $\sqrt{n} \approx 14$ and thus we can set $u_A = u$. We use $\Vert \vert  \Lbar \vert \vert \Lbar^{-1} \vert \Vert + \Vert \vert  \Ubar^{-1} \vert \vert \Ubar \vert \Vert $ to approximate $\Vert E_R \Vert / \Vert M_R^{-1} \Vert$ in Table~\ref{tab:left-right-precond_condition_numbers}. As in the split-preconditioning case, the values are larger than $u_R^{-1}$ for half precision.

Experiments with $u=u_A=u_L$ set to double (Table~\ref{tab:results_uL_double_left_precond}) show that the bounds for $\psi_L$ in Table~\ref{tab:left-right-precond_condition_numbers} largely overestimate the error in applying the preconditioner and even though the bounds for $\psi_A$ are quite close to the obtained values, we still obtain $\mathcal{O}(u)$ relative backward error for the unpreconditioned system. Note that the $\zeta$ values are increased by $\psi_A$. If we keep $u=u_A$ set to double and set $u_L$ to single, the relative backward error for the unpreconditioned system reaches $\mathcal{O}(u_L)$ and $u_L \psi_L$ determines $\zeta$ (Table~\ref{tab:results_uL_single}). This agrees with the split-reconditioning results. 

We keep $u=u_A$ set to double for right-preconditioning experiments (Tables~\ref{tab:results_uR_double_right_precond} and \ref{tab:results_uR_single_right_precond}). Note that the relative backward error reaches $\mathcal{O}(u)$ with $u_R$ set to both double and single, except for $\kappa(A)= 10^{10}$ with $u_R$ set to single. The number of iterations when $u_R$ is set to single grows for highly ill-conditioned problems. Though the term $\Vert \Zbar_k \Vert \Vert M_R(\xbar_k - \xbar_0) \Vert$ grows as in the split-preconditioning case, it is balanced by $\Vert \xbar_k \Vert$ here. Note that in these experiments the bound in \cite[eq. (3.22) ]{arioli2009using} is applicable for $c\geq 6$. Here $\Gamma = \Vert \vert L \vert \vert U \vert \Vert / \Vert A \Vert = \mathcal{O}(10)$ and the term $\sqrt{u} \Gamma \Vert A \Vert \Vert \Zbar_k \Vert$ is thus small ensuring a tight bound for the backward error. 

\begin{figure}
    \centering
\begin{subfigure}[t]{0.45\linewidth}
  \centering
 \includegraphics[width=\linewidth]{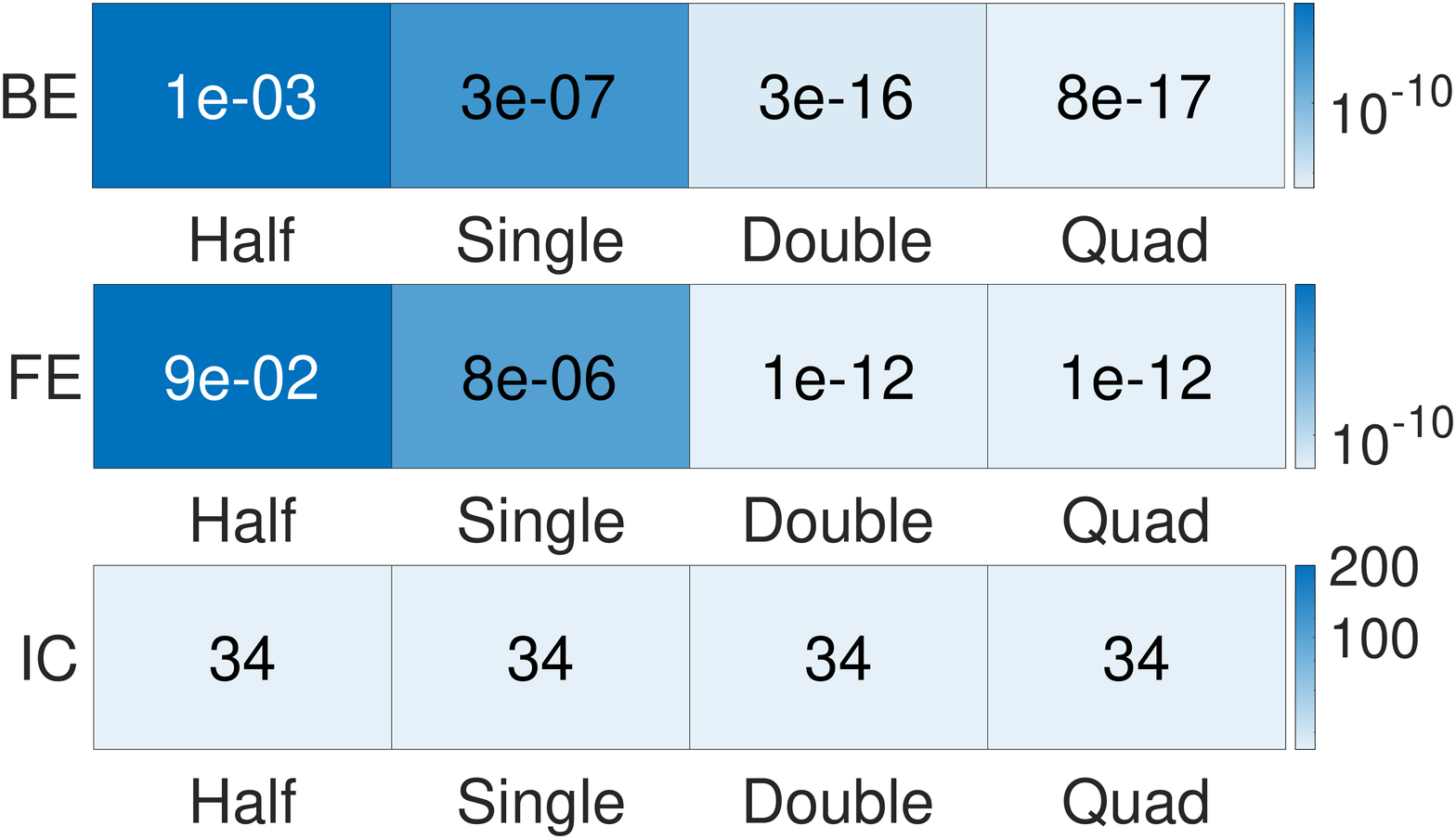}
  \caption{Left-preconditioning}
\end{subfigure}
\begin{subfigure}[t]{0.45\linewidth}
  \centering
 \includegraphics[width=\linewidth]{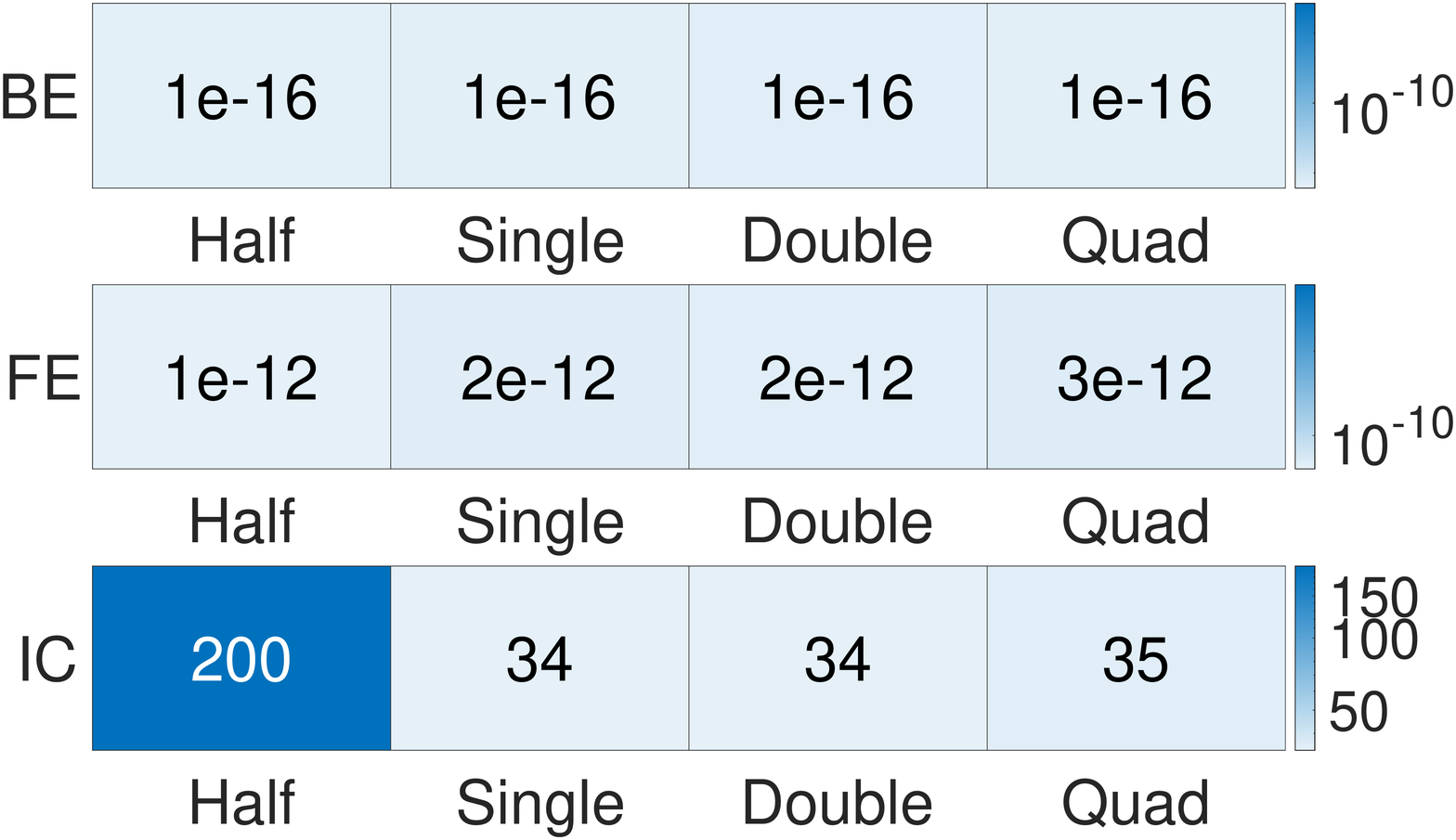}
  \caption{Right-preconditioning}
\end{subfigure}
    \caption{As in Figure~\ref{fig:dense_c5}, but for full left- and right-preconditioning. The left panel shows results for different choices of $u_L$ and the right panel shows results for different choices of $u_R$. IC is the iteration count.}
    \label{fig:dense_c5_full_left_right}
\end{figure}

The results for $c=5$ for all choices of $u_L$ and $u_R$ are presented in Figure~\ref{fig:dense_c5_full_left_right}. 
Comparing these with the skew diagonals (/) of the respective heatmaps in Figure~\ref{fig:dense_c5}, that is, when the preconditioners are applied in the same precision, shows that applying the full preconditioner on the left gives larger relative backward and consequentially forward errors when the preconditioner is applied in low precision. This is not the case when the full preconditioner is applied on the right. It is curious that the iteration count with left-, right- and split-preconditioning is essentially different only when the preconditioner is applied in half precision. The results for different $c$ values are similar. 

\begin{table}[]
    \caption{Synthetic problems. Condition numbers of the preconditioned coefficient matrices with full left- and right-preconditioning, and bounds for $\psi_A$ and $\psi_L$ (left-preconditioning only). The preconditioners are computed in precision accurate to four decimal digits for $c<6$ and in single precision for $c\geq 6$.}
    \centering
    \begin{tabular}{c|c|c|c|c|c} 
    & \multicolumn{3}{c|}{left-preconditioning} &\multicolumn{2}{c}{right-preconditioning} \\
c & $\kappa(\Ubar^{-1} \Lbar^{-1} A)$	& $\psi_A$ bound & $\psi_L$ bound & $\kappa(  A \Ubar^{-1} \Lbar^{-1})$ & $\Vert E_R \Vert / \Vert M_R^{-1} \Vert$ approx. \\	
\hline
1 & $1.01$ & $2.76 \times 10^2 $ & $3.79 \times 10^5 $ & $1.01$ & $2.74 \times 10^3 $\\ 
2 & $1.02$ & $1.40 \times 10^3 $ & $6.45 \times 10^5 $ & $1.02$ & $2.84 \times 10^3 $\\
3 & $1.09 $ & $9.17 \times 10^3 $ & $2.35 \times 10^6$ & $1.09$ & $3.24 \times 10^3 $\\
4 & $1.75$ &  $5.53 \times 10^4 $ &  $1.51 \times 10^7$ & $1.74$ & $3.27 \times 10^3 $\\
5 & $3.09 \times 10^1$ &  $9.54 \times 10^4 $ &  $1.93 \times 10^7$ & $3.09 \times 10^1$ & $3.30 \times 10^3 $\\
\hline
6 & $1.03 $ & $4.90 \times 10^6 $ & $6.04 \times 10^8$ & $1.03$ & $3.39 \times 10^3 $\\
7 & $1.31$ & $3.79 \times 10^7$ & $3.12 \times 10^9$ & $1.30$ & $3.66 \times 10^3 $\\
8 & $7.79$ & $1.39 \times 10^8$ & $1.93 \times 10^{10} $ & 7.56 & $3.96 \times 10^3 $\\
9 & $4.75 \times 10^3$ & $8.72 \times 10^7$ & $7.14 \times 10^9$ & $1.96 \times 10^3$ & $4.08 \times 10^3 $\\
10 & $2.49 \times 10^5$ & $1.01 \times 10^8$ & $3.75 \times 10^9 $ & $4.03 \times 10^4 $ & $4.08 \times 10^3$
\end{tabular}
\label{tab:left-right-precond_condition_numbers}
\end{table}

\begin{table}[]
    \caption{Synthetic problems with full left-preconditioning and $u_L$ set to double. IC denotes the iteration count, BE is the relative backward error and FE is the relative forward error.}
    \centering
     \scalebox{0.7}{
    \begin{tabular}{c|c|c|c|c|c|c}
c & IC & BE &	FE	&	$\zeta$ & $\psi_A$	& $\psi_L$  \\
\hline
1 & 6  & $1.60 \times 10^{-15} $ & $6.77 \times 10^{-15} $ & $5.21 \times 10^{-15} $  & 6.41 & $4.68 \times 10^1$  \\
2 & 7 & $9.17 \times 10^{-16}$ & $8.39 \times 10^{-15} $ & $8.42 \times 10^{-15}$ & $2.99 \times 10^1 $ & $4.64 \times 10^1 $\\
3 & 9 & $4.56 \times 10^{-16} $ & $2.31 \times 10^{-14}$ & $2.66 \times 10^{-14} $ & $1.95 \times 10^2$ & $4.41 \times 10^1$ \\
4 & 15 & $3.53 \times 10^{-16} $ & $1.61 \times 10^{-13}$ & $1.55 \times 10^{-13} $ & $1.36 \times 10^3 $ & $4.23 \times 10^1$\\
5 & 34 & $3.33 \times 10^{-16} $ & $1.06 \times 10^{-12} $ & $2.63 \times 10^{-13} $ & $2.35 \times 10^3 $ & $1.35 \times 10^1 $ \\
\hline
6 & 7 & $1.31 \times 10^{-16}$ & $1.14 \times 10^{-11}$ & $1.13 \times 10^{-11} $ & $1.01 \times 10^5$ & $4.26 \times 10^1$ \\
7 & 11 & $1.17 \times 10^{-16}$ & $7.71 \times 10^{-11}$ & $1.23 \times 10^{-10}$ & $1.10 \times 10^6 $ & $5.43 \times 10^1$ \\
8 & 21 & $1.12 \times 10^{-16} $ & $5.25 \times 10^{-10}$ & $4.03 \times 10^{-10}$ & $3.63 \times 10^6$ & $4.50 \times 10^1$\\
9 & 47 & $3.90 \times 10^{-16}$ & $8.10 \times 10^{-9}$ & $3.27 \times 10^{-10} $ & $2.94 \times 10^6$ & $8.02 \times 10^1 $ \\
10 & 124 & $6.69 \times 10^{-16} $ & $1.42 \times 10^{-7} $ & $5.78 \times 10^{-10}$ & $5.20 \times 10^6$ & $9.46 \times 10^3 $
    \end{tabular}
    }
    \label{tab:results_uL_double_left_precond}
\end{table}

\begin{table}[]
    \caption{As in Table~\ref{tab:results_uL_double_left_precond}, but for $u_L$ set to single. IC denotes the iteration count, BE is the relative backward error and FE is the relative forward error.}
    \centering
     \scalebox{0.7}{
    \begin{tabular}{c|c|c|c|c|c|c}
c & IC & BE &	FE	&	$\zeta$ & $\psi_A$	& $\psi_L$  \\
\hline
1 & 6 & $9.29 \times 10^{-7}$ & $3.85 \times 10^{-6}$ & $2.32 \times 10^{-6} $ & 6.48 & $4.52 \times 10^1$\\
2 & 7 & $6.07 \times 10^{-7} $ & $4.25 \times 10^{-6} $ & $2.56 \times 10^{-6}$ & $2.86 \times 10^1 $ & $4.38 \times 10^1$ \\
3 & 9 & $2.60 \times 10^{-7}$  & $ 2.76 \times 10^{-6}$ & $2.55 \times 10^{-6}$ & $2.07 \times 10^2 $ & $4.29 \times 10^1$ \\
4 & 15 & $1.72 \times 10^{-7} $ & $4.97 \times 10^{-6}$ & $2.15 \times 10^{-6} $ & $1.32 \times 10^{3} $ & $3.61 \times 10^1 $ \\
5 & 34 & $2.78 \times 10^{-7} $ & $7.99 \times 10^{-6}$ & $8.52 \times 10^{-7}$ & $3.08 \times 10^3 $ & $1.43 \times 10^1 $ \\
\hline
6 & 7 & $1.00 \times 10^{-7} $ & $3.99 \times 10^{-6}$ & $2.84 \times 10^{-6}$ & $1.06 \times 10^5$ & $4.76 \times 10^1$ \\
7 & 11 & $4.41 \times 10^{-8}$ & $4.41 \times 10^{-8}$ & $3.52 \times 10^{-6}$ & $9.85 \times 10^5$ & $5.90 \times 10^1$ \\
8 & 21 & $5.01 \times 10^{-8} $ & $4.46 \times 10^{-6}$ & $ 2.67 \times 10^{-6}$ & $3.37 \times 10^6$ & $4.48 \times 10^1$ \\
9 & 47 & $1.55 \times 10^{-7}$ & $3.23 \times 10^{-5}$ & $ 3.98 \times 10^{-6} $ & $3.99 \times 10^6$ & $6.67 \times 10^1$ \\
10 & 114 & $2.56 \times 10^{-7}$ & $2.74 \times 10^{-2}$ & $4.55 \times 10^{-3}$ & $5.78 \times 10^6$ & $7.63 \times 10^4$
    \end{tabular}
    }
    \label{tab:results_uL_single_left_precond}
\end{table}

\begin{table}[]
    \caption{Synthetic problems with full right-preconditioning and $u_R$ set to double. IC denotes the iteration count, BE is the relative backward error and FE is the relative forward error.}
    \centering
     \scalebox{0.7}{
    \begin{tabular}{c|c|c|c|c|c|c|c}
c	& IC & BE &	FE	&	$\zeta$&	$\Vert \Zbar_k \Vert \Vert M_R(\xbar_k - \xbar_0) \Vert$	& $\psi_A$	& $\rho$  \\
\hline
1 & 6 & $1.18 \times 10^{-16}$ & $6.83 \times 10^{-16} $ & $4.31 \times 10^{-16} $ & $4.24 \times 10^1$ & 1.41 & $5.96 \times 10^{-16}$ \\
2 & 7 & $9.15 \times 10^{-17} $ & $2.73 \times 10^{-15}$ & $5.10 \times 10^{-16}$ & $3.34 \times 10^2$ & $9.58 \times 10^{-1}$ & $ 4.69 \times 10^{-15} $ \\
3 & 9 & $8.00 \times 10^{-17}$ & $1.82 \times 10^{-14}$ & $6.39 \times 10^{-16}$ & $3.47 \times 10^3$ & $8.48 \times 10^{-1} $ & $4.88 \times 10^{-14}$ \\
4 & 15 & $8.13 \times 10^{-17}$ & $1.50 \times 10^{-13}$ & $1.50 \times 10^{-13} $ & $3.26 \times 10^4$ & $6.53 \times 10^{-1} $ & $4.60 \times 10^{-13} $ \\
5 & 34 & $1.39 \times 10^{-16}$ & $1.65 \times 10^{-12}$ & $1.62 \times 10^{-15} $ & $7.79 \times 10^5 $ & $5.95 \times 10^{-1} $ & $1.01 \times 10^{-11} $ \\
\hline
6 & 7 & $5.08 \times 10^{-17} $ & $1.04 \times 10^{-11}$ & $5.77 \times 10^{-16} $ & $2.49 \times 10^6 $ & $6.09 \times 10^{-1}$ & $3.50 \times 10^{-11} $ \\
7 & 11 & $5.93 \times 10^{-17}$ & $7.32 \times 10^{-11}$ & $4.66 \times 10^{-16}$ & $1.92 \times 10^7 $ & $5.48 \times 10^{-1}$ & $2.69 \times 10^{-10}$ \\
8 & 21 & $4.87 \times 10^{-17}$ & $7.15 \times 10^{-10}$ & $9.40 \times 10^{-16}$ & $3.76 \times 10^8 $ & $4.71 \times 10^{-1} $ & $5.25 \times 10^{-9} $ \\
9 & 61 & $3.13 \times 10^{-16}$ & $4.04 \times 10^{-8} $ & $9.87 \times 10^{-15} $ & $3.77 \times 10^{10} $ & $4.32 \times 10^{-1}$ & $3.06 \times 10^{-7}$ \\
10 & 200 & $8.15 \times 10^{-16} $ & $1.44 \times 10^{-6} $ & $4.71 \times 10^{-14}$ & $1.53 \times 10^{12} $ & $5.77 \times 10^{-1} $ & $1.62 \times 10^{-6}$
    \end{tabular}
    }
    \label{tab:results_uR_double_right_precond}
\end{table}

\begin{table}[]
    \caption{As in Table~\ref{tab:results_uR_double_right_precond}, but for $u_R$ set to single. IC denotes the iteration count, BE is the relative backward error and FE is the relative forward error.}
    \centering
     \scalebox{0.7}{
    \begin{tabular}{c|c|c|c|c|c|c|c}
c	& IC & BE &	FE	&	$\zeta$&	$\Vert \Zbar_k \Vert \Vert M_R(\xbar_k - \xbar_0) \Vert$	& $\psi_A$	& $\rho$  \\
\hline
1 & 6 & $1.21 \times 10^{-16} $ & $ 6.77 \times 10^{-16} $ & $4.18 \times 10^{-16} $ & $4.25 \times 10^1$ & $1.34$ & $1.09 \times 10^{-7}$\\
2 & 7 & $9.06 \times 10^{-17}$ & $3.21 \times 10^{-15} $ & $4.88 \times 10^{-16}$ & $3.34 \times 10^2 $ & $8.72 \times 10^{-1}$ & $1.94 \times 10^{-7}$ \\
3 & 9 & $7.87 \times 10^{-17}$ & $1.88 \times 10^{-14}$ & $6.15 \times 10^{-16}$ & $3.48 \times 10^3 $ & $7.73 \times 10^{-1}$ & $7.15 \times 10^{-7}$\\
4 & 15 & $7.62 \times 10^{-17} $ & $1.25 \times 10^{-13} $ & $6.32 \times 10^{-16}$ & $3.26 \times 10^4$ & $6.72 \times 10^{-1}$ & $4.98 \times 10^{-6}$ \\
5 & 34 & $1.15 \times 10^{-16}$ & $1.91 \times 10^{-12}$ & $1.63 \times 10^{-15}$ & $7.79 \times 10^5$ & $6.05 \times 10^{-1}$ & $2.39 \times 10^{-5}$ \\
\hline
6 & 7 & $5.31 \times 10^{-17}$ & $9.20 \times 10^{-12} $ & $9.20 \times 10^{-12} $ & $1.25 \times 10^6 $ & $4.42 \times 10^{-1} $ & $1.48 \times 10^{-4}$\\
7 & 11 & $4.63 \times 10^{-17}$ & $6.76 \times 10^{-11}$ & $3.39 \times 10^{-16}$ & $1.42 \times 10^7$ & $5.21 \times 10^{-1}$ & $7.40 \times 10^{-4}$ \\
8 & 27 & $6.27 \times 10^{-17}$ & $5.69 \times 10^{-10} $ & $5.82 \times 10^{-16} $ & $2.35 \times 10^8 $ & $4.58 \times 10^{-1}$ & $4.58 \times 10^{-1}$ \\
9 & 200 & $1.97 \times 10^{-16} $ & $2.25 \times 10^{-8}$ & $1.86 \times 10^{-14}$ & $6.94 \times 10^{10}$ & $4.71 \times 10^{-1} $ & $1.77 \times 10^{-1} $ \\
10 & 200 & $1.41 \times 10^{-15}$ & $2.29 \times 10^{-6}$ & $3.98 \times 10^{-14} $ & $1.32 \times 10^{12} $ & $5.43 \times 10^{-1} $ & $2.22 \times 10^{-1} $
    \end{tabular}
    }
    \label{tab:results_uR_single_right_precond}
\end{table}

\subsection{Numerical example: application problems with split-preconditioning}\label{sec:numerics_sparse}

We perform experiments with some ill-conditioned problems from the SuiteSparse collection \cite{SuiteSparse}; see Table~\ref{tab:sparsesuit_problems}. We generate the right hand side vector $b$ in the same way as for the synthetic problems. The preconditioner is computed as in the previous section, but in single precision, and the matrix $A$ is converted to a full matrix due to the lacking single precision sparse matrix-vector product functionality in MATLAB. We report results for split-preconditioning only; these illustrate the trends in left- and right-preconditioning cases well. Note that the split-preconditioner reduces the condition number to the theoretical minimum or close to it, except for the problem with the highest $\kappa(A)$. However $\kappa(\Atilde)$ and $\kappa(M_R)$ are of the order of $\kappa(A)$ and thus we expect to see its effect on the forward error. We set $u_A$ and $u$ to double based on the $\psi_A$ bound; $u_L$ and $u_R$ are set as in the previous sections. Unpreconditioned FGMRES does not converge in $n$ iterations for any of the problems. 

\begin{figure}
    \centering
\begin{subfigure}[t]{0.45\linewidth}
  \centering
 \includegraphics[width=\linewidth]{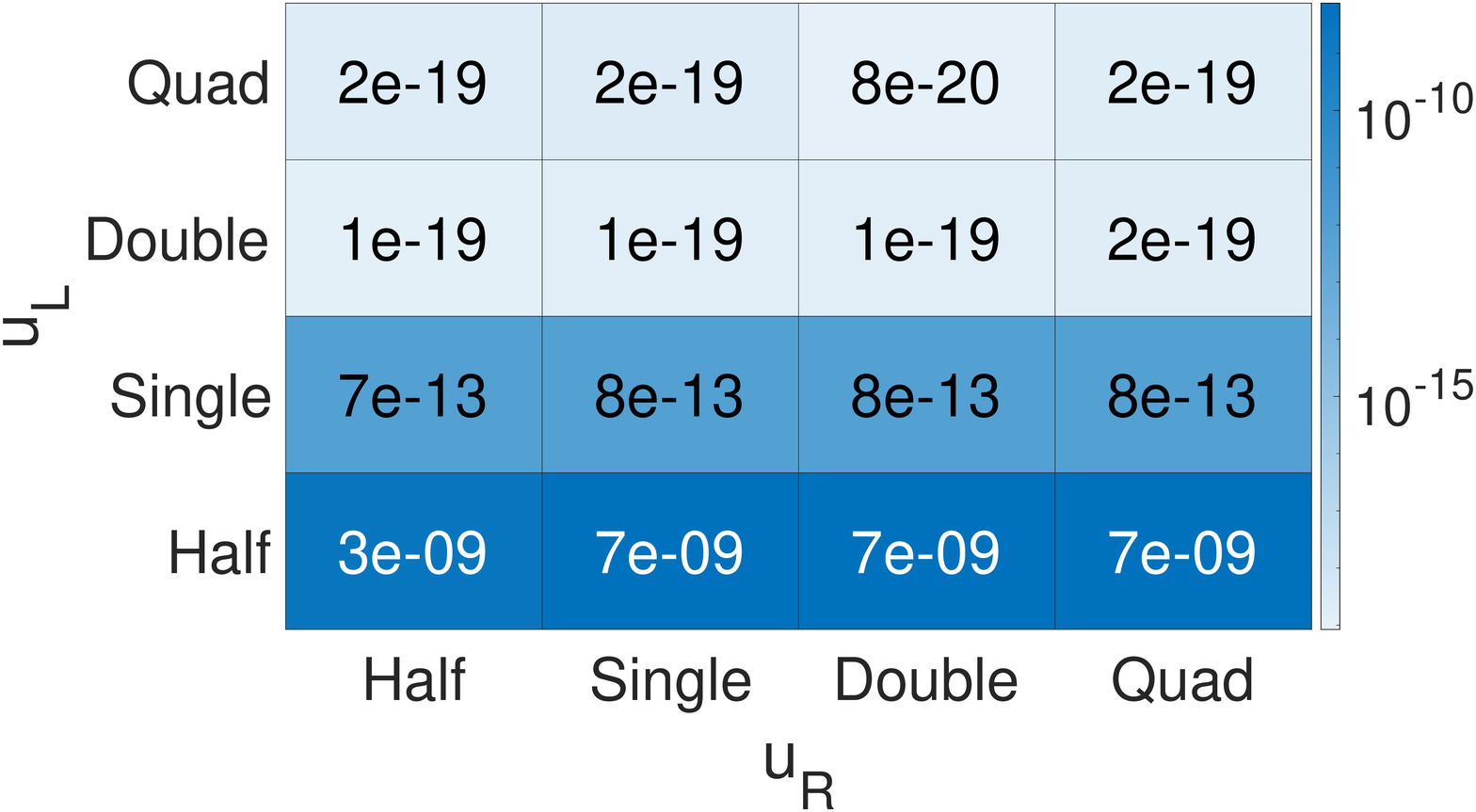}
  \caption{BE, rajat14}
\end{subfigure}
\begin{subfigure}[t]{0.45\linewidth}
  \centering
 \includegraphics[width=\linewidth]{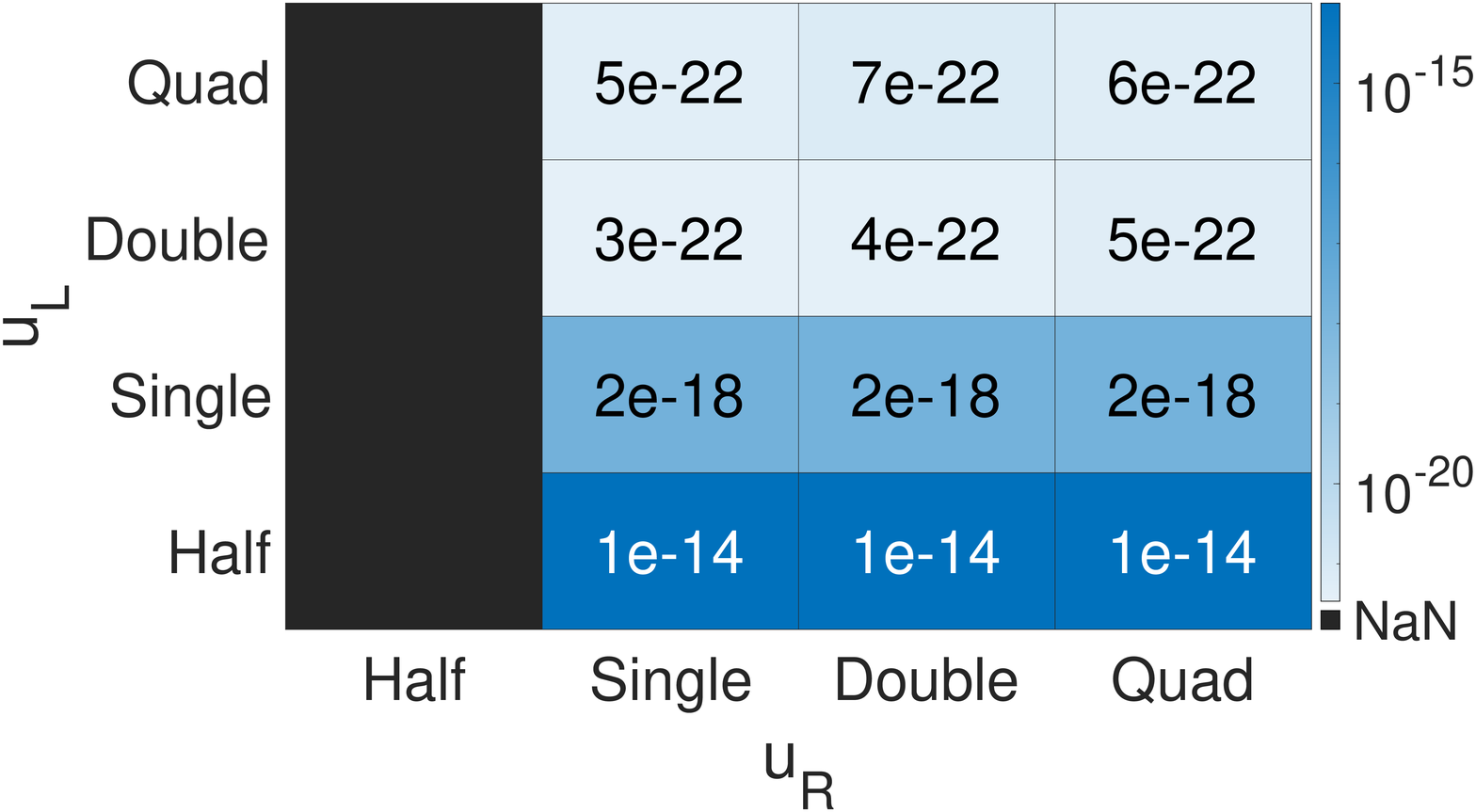}
  \caption{BE, arc130}
\end{subfigure}
\begin{subfigure}[t]{0.45\linewidth}
  \centering
 \includegraphics[width=\linewidth]{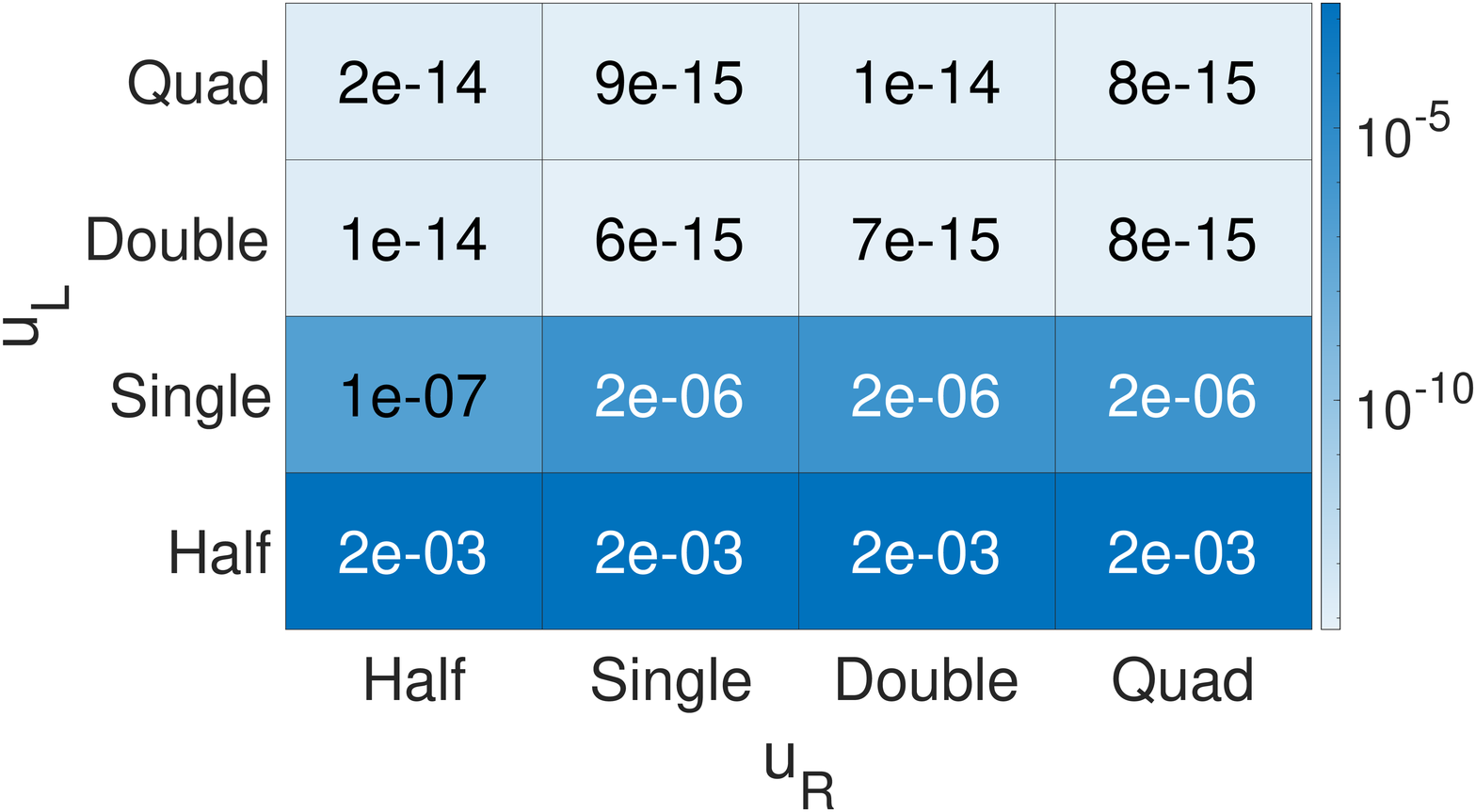}
  \caption{FE, rajat14}
\end{subfigure}
\begin{subfigure}[t]{0.45\linewidth}
  \centering
 \includegraphics[width=\linewidth]{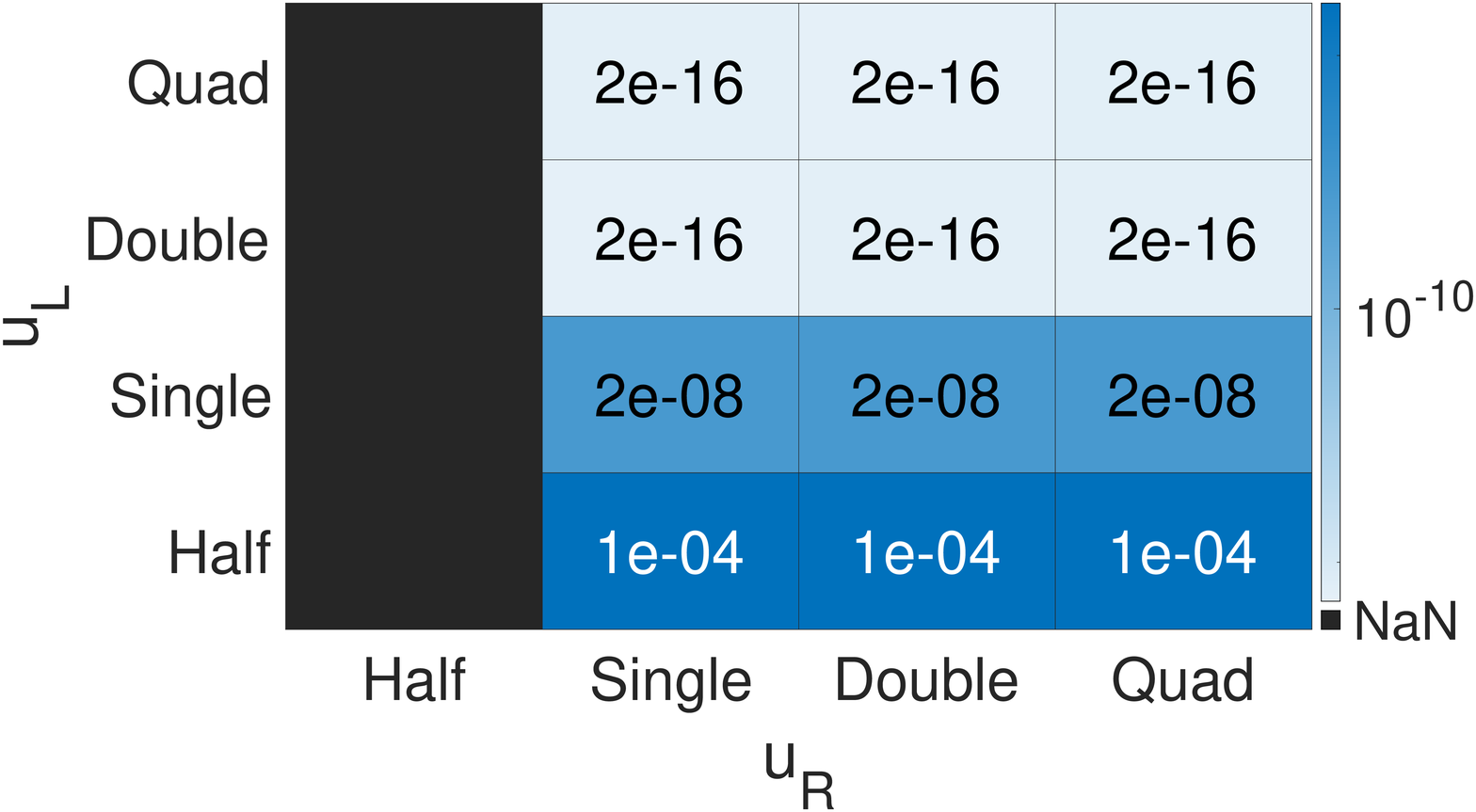}
  \caption{FE, arc130}
\end{subfigure}
\begin{subfigure}[t]{0.45\linewidth}
  \centering
 \includegraphics[width=\linewidth]{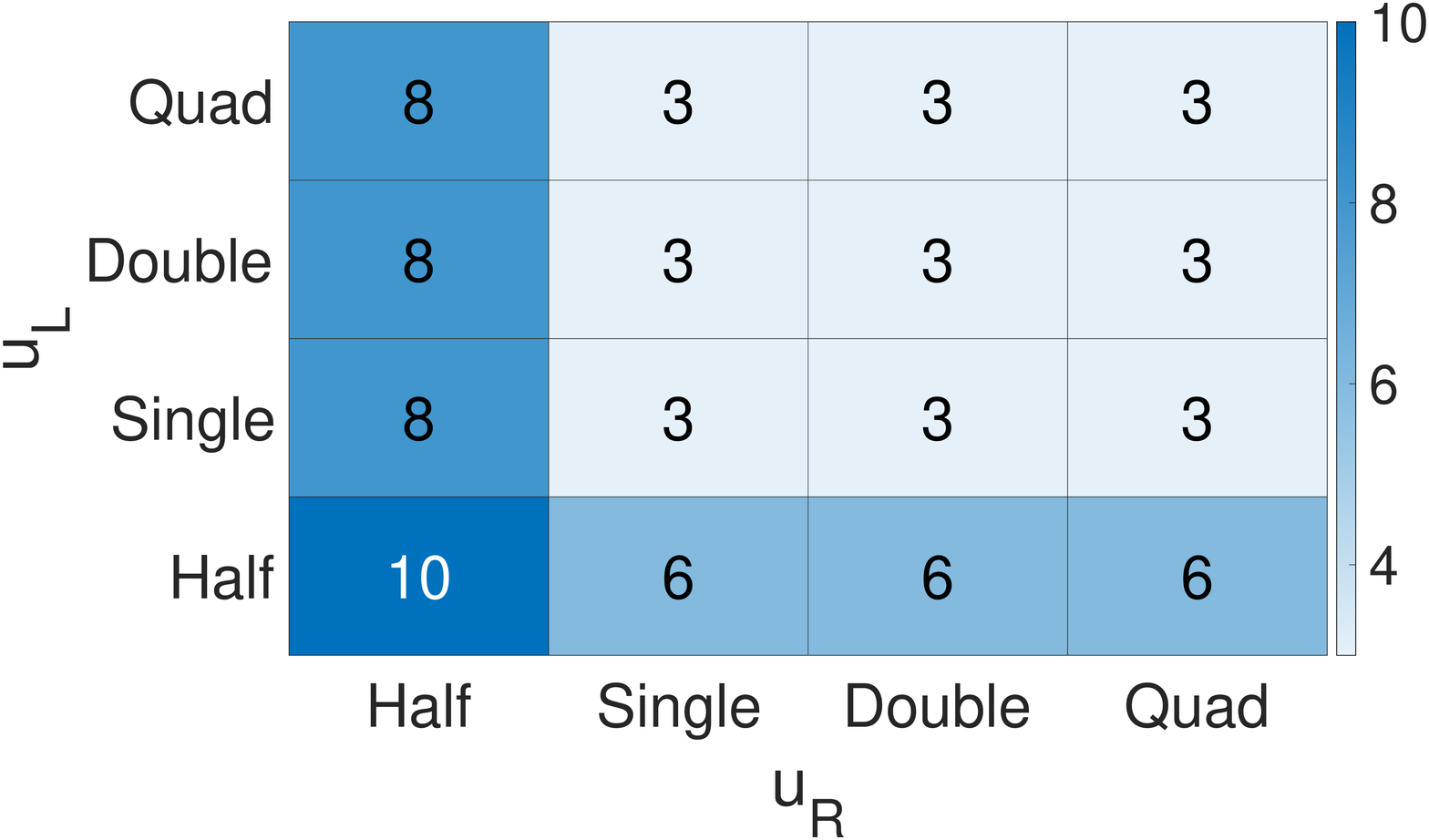}
  \caption{iterations, rajat14}
\end{subfigure}
\begin{subfigure}[t]{0.45\linewidth}
  \centering
 \includegraphics[width=\linewidth]{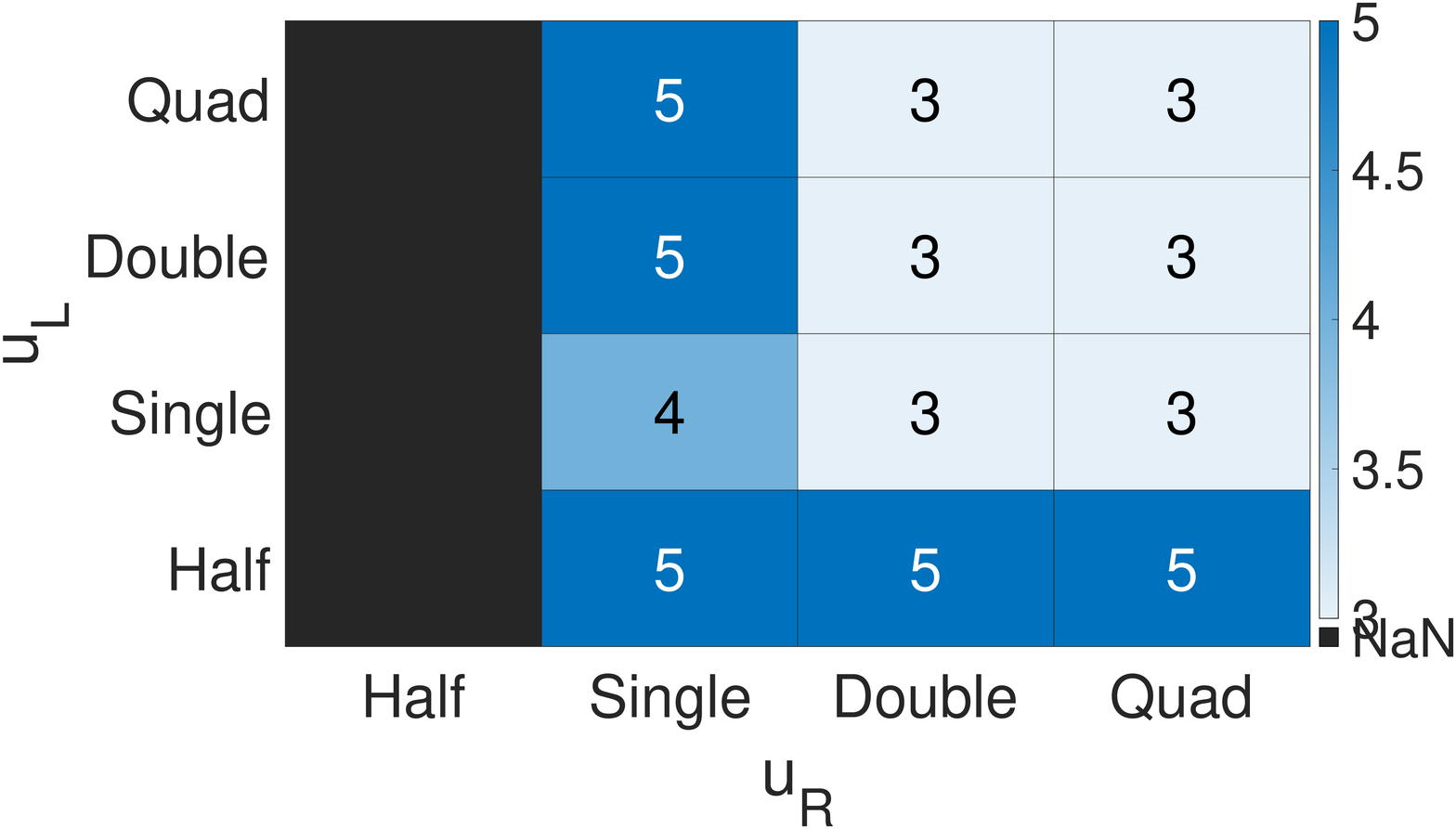}
  \caption{iterations, arc130}
\end{subfigure}
    \caption{SuiteSparse problems rajat14 and arc130. BE is the relative backward error and FE is the relative forward error.}
    \label{fig:sparse_rajat14_arc130}
\end{figure}

\begin{figure}
    \centering
\begin{subfigure}[t]{0.45\linewidth}
  \centering
 \includegraphics[width=\linewidth]{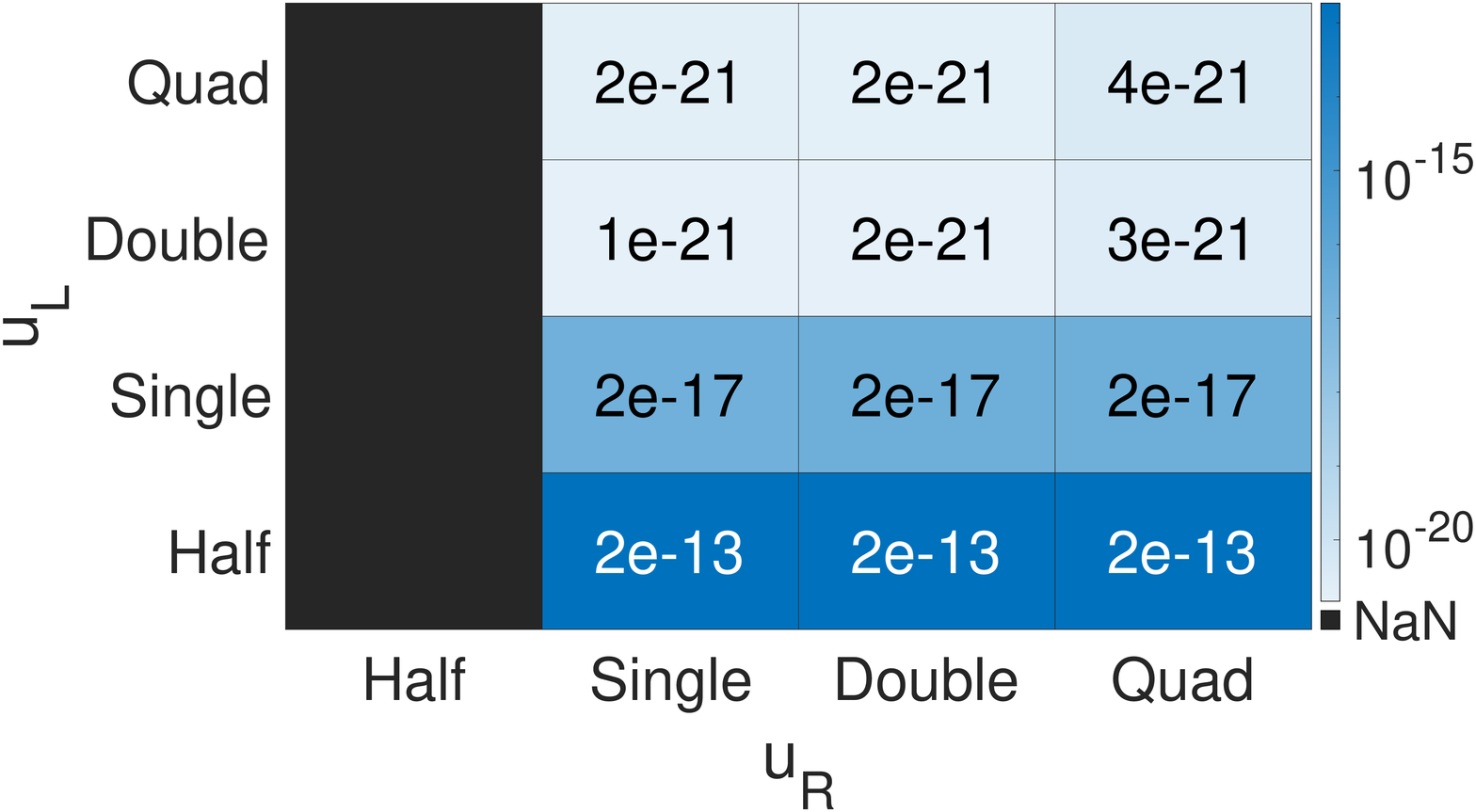}
  \caption{BE, west0132}
\end{subfigure}
\begin{subfigure}[t]{0.45\linewidth}
  \centering
 \includegraphics[width=\linewidth]{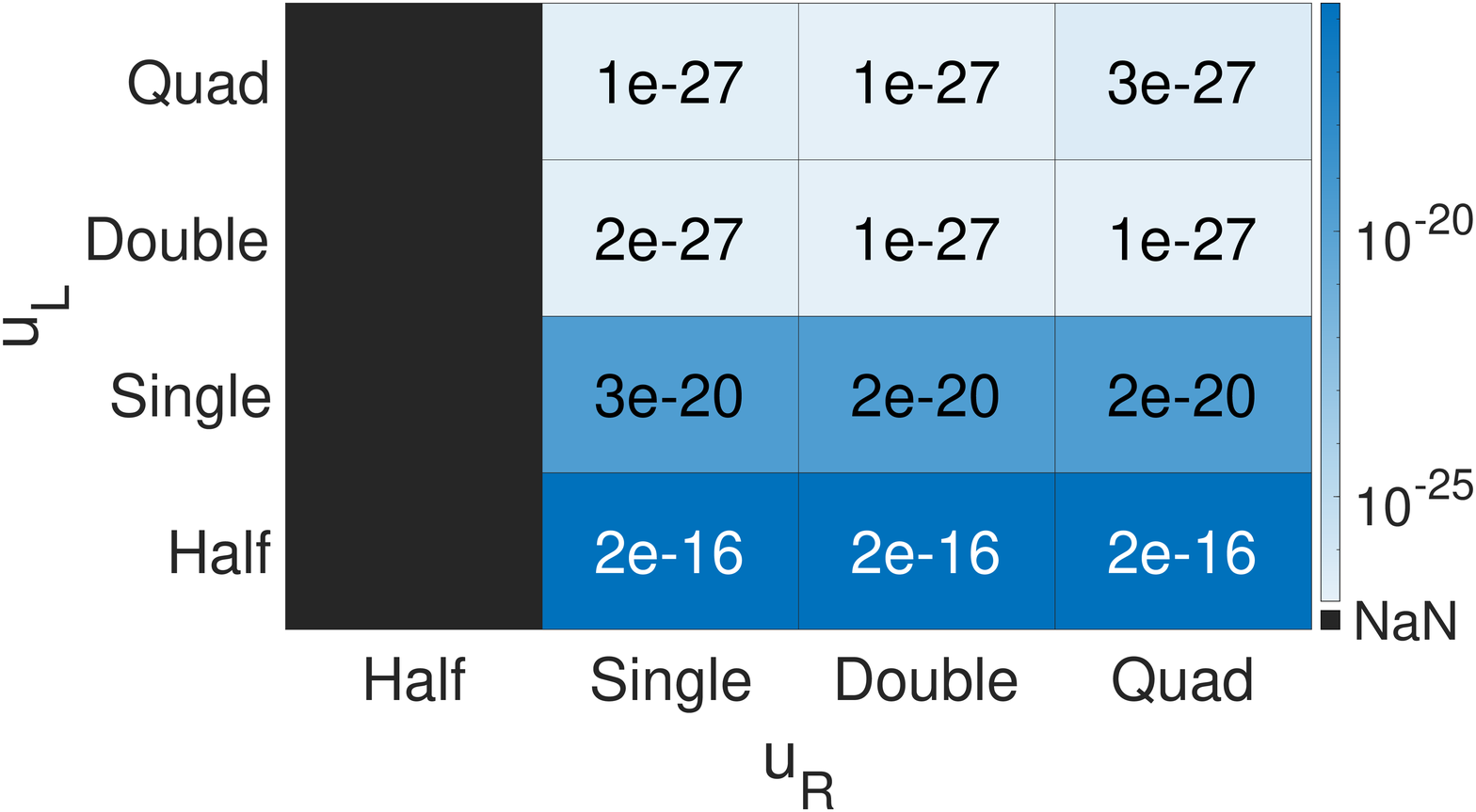}
  \caption{BE, fs\_183\_3}
\end{subfigure}
\begin{subfigure}[t]{0.45\linewidth}
  \centering
 \includegraphics[width=\linewidth]{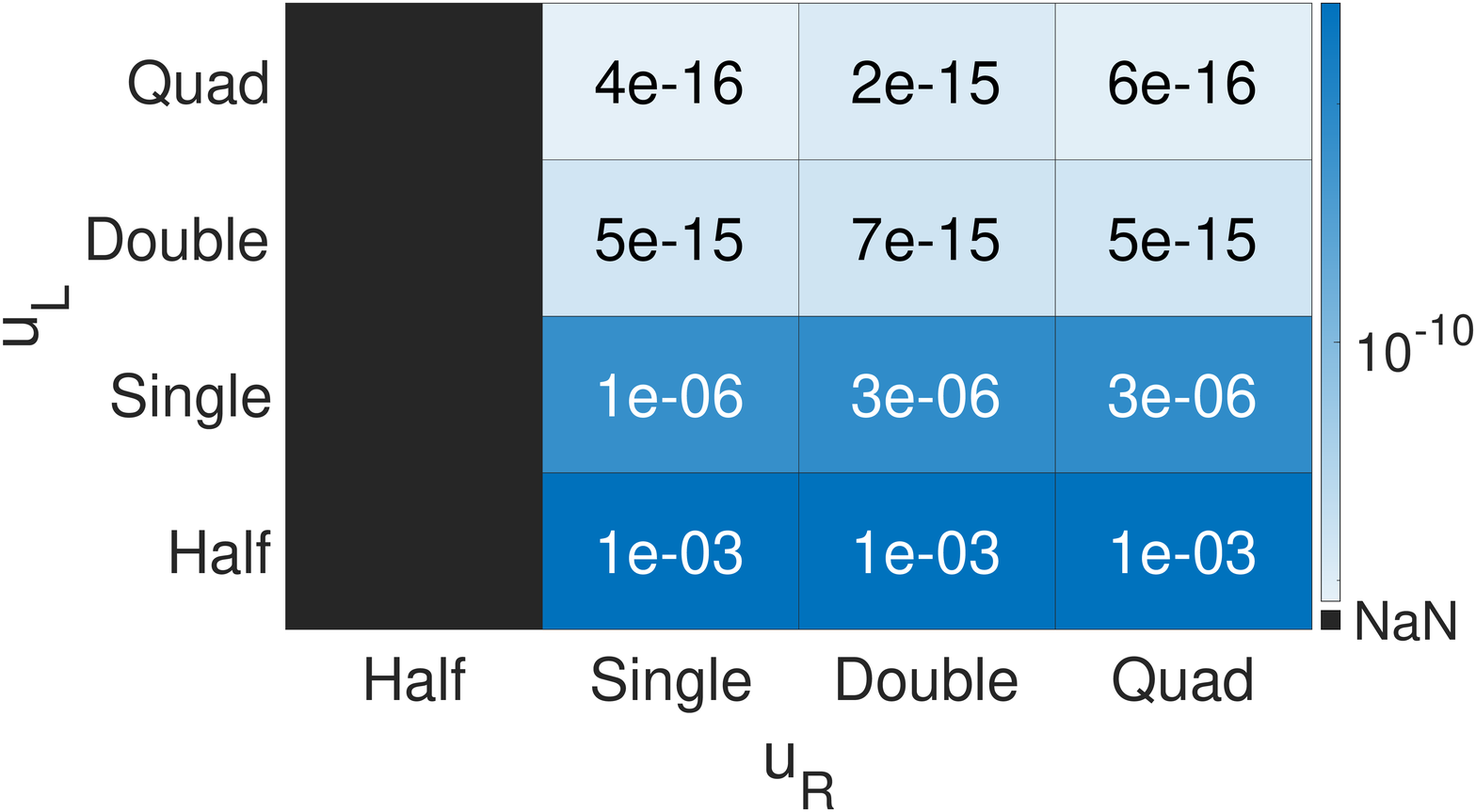}
  \caption{FE, west0132}
\end{subfigure}
\begin{subfigure}[t]{0.45\linewidth}
  \centering
 \includegraphics[width=\linewidth]{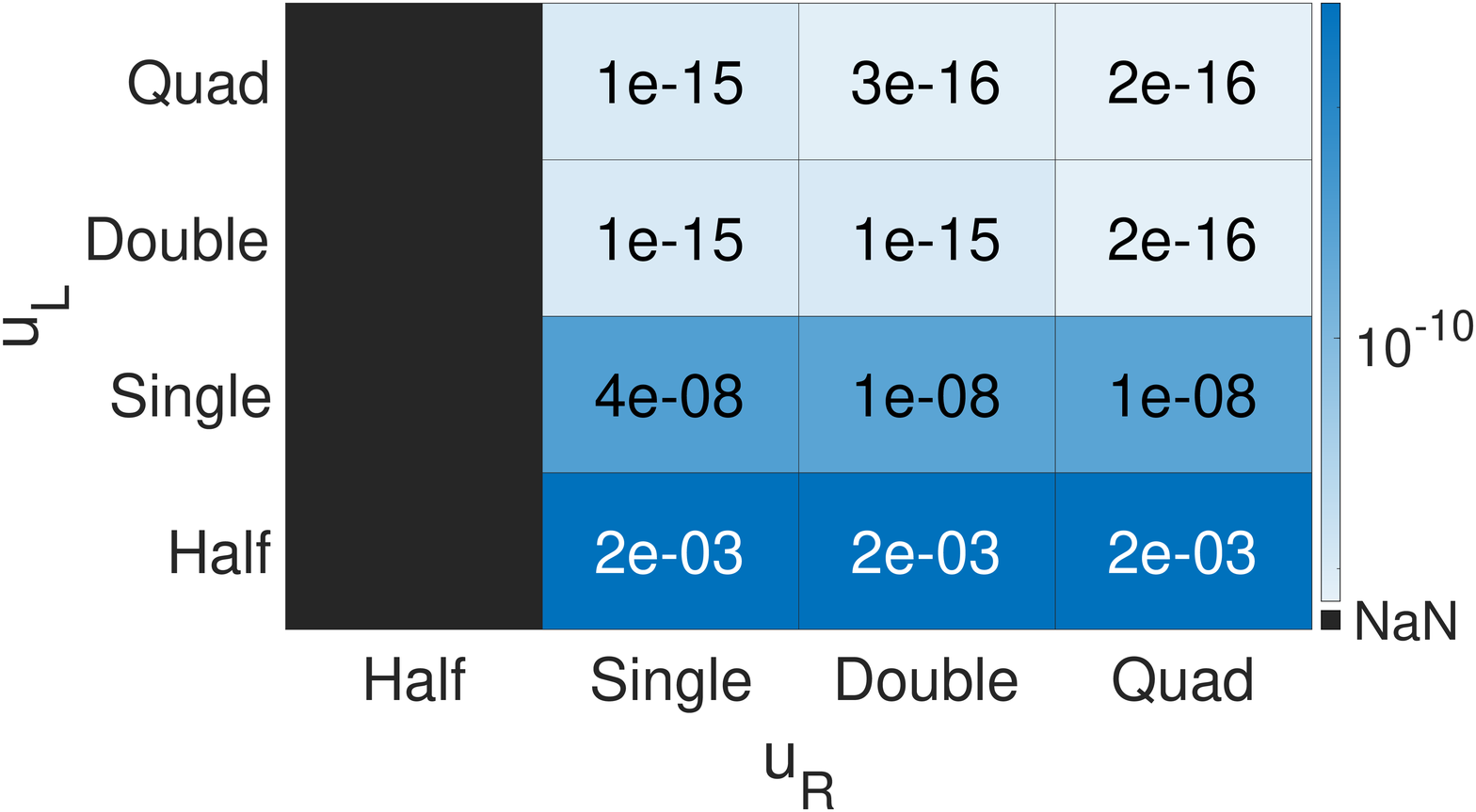}
  \caption{FE, fs\_183\_3}
\end{subfigure}
\begin{subfigure}[t]{0.45\linewidth}
  \centering
 \includegraphics[width=\linewidth]{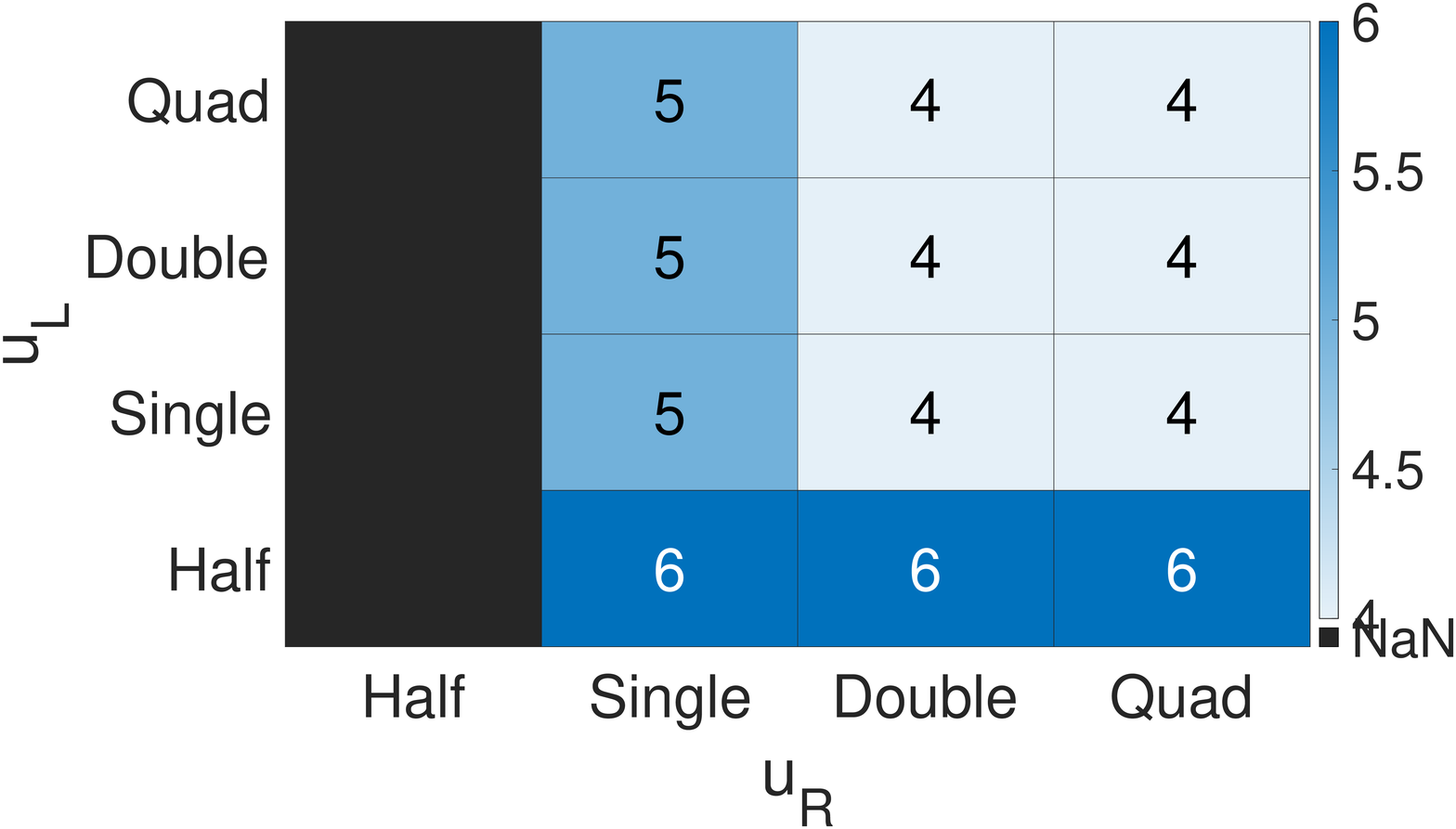}
  \caption{iterations, west0132}
\end{subfigure}
\begin{subfigure}[t]{0.45\linewidth}
  \centering
 \includegraphics[width=\linewidth]{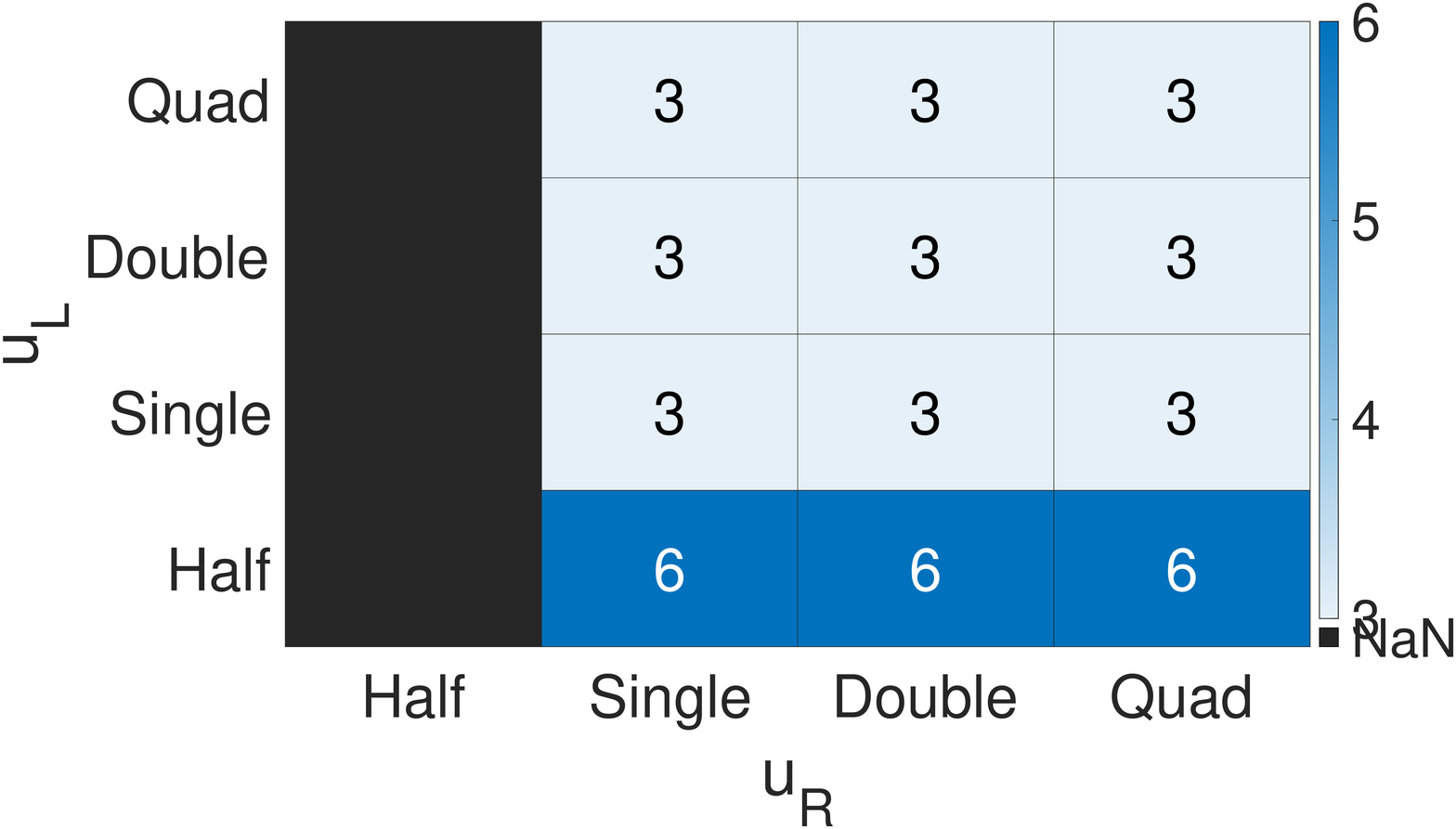}
  \caption{iterations, fs\_183\_3}
\end{subfigure}
    \caption{SuiteSparse problems west0132 and fs\_183\_3. BE is the relative backward error and FE is the relative forward error.}
    \label{fig:sparse_wets0132_fs1833}
\end{figure}

Approximations of $\Vert E_R \Vert / \Vert M_R^{-1} \Vert$ indicate that the backward error may be affected if we apply $M_R$ in half precision for \textit{arc130} and \textit{west0132}. However, we cannot test this as for all problems except \textit{rajat14}, $M_R$ is singular with respect to $u_R$ set to half. This may be amended by using scaling strategies when computing the preconditioner; see, for example, \cite{HighamSqueeze}. We observe similar tendencies (Figures~\ref{fig:sparse_rajat14_arc130} and \ref{fig:sparse_wets0132_fs1833}) as for the dense problems, however here we can achieve smaller backward error and for \textit{fs\_183\_3} the backward error is $\mathcal{O}(u)$ even with $u_L$ set to half. Note that for sparse problems setting $u_L$ to a low precision results in iteration-wise slower convergence. The term $\Vert \Zbar_k \Vert \Vert M_R(\xbar_k - \xbar_0) \Vert$ grows as for the dense problems, but is balanced by $\Vert \xbar_k \Vert$; see Tables~\ref{tab:sparsesuit_uL_single} and \ref{tab:sparsesuit_uR_single}.

In all of our sparse and dense examples, $\psi_A/\psi_L$ does not become large enough to allow setting $u_L > u_A$ without it affecting the backward error. However for \textit{fs\_183\_3} both $\psi_A$ and $\psi_L$ are small enough that we can set $u_A$ and $u_L$ to single and expect $\mathcal{O}(u)$ backward error. Numerical experiments confirm this even though the backward and forward errors become slightly larger compared to $u_A$ set to double (not shown).

\begin{table}[]
\caption{As in Table~\ref{tab:cond_no_different_c}, but for SuiteSparse problems.}
    \centering
       \scalebox{0.8}{
    \begin{tabular}{c|c|c|c|c|c|c|c|c}
 problem	& $n$ & $\kappa(A)$ &$\kappa(\Atilde)$ &	$\kappa(\Ahat)$&	$\kappa(M_R)$	&$\kappa(M_L)$	& $\psi_A$ bound &$\frac{\Vert E_R \Vert }{ \Vert M_R^{-1} \Vert} $  \\	
& & & & & & & & approx. \\
 \hline
rajat14	& 180 &$3.22 \times 10^8$	&$1.44 \times 10^8$& 1.01	& $1.44 \times 10^8$& 	$9.72 \times 10^1$&	1.13 & 33 \\
arc130 & 130 &	$6.05 \times 10^{10}$&	$6.05 \times 10^{10}$&	1.00	& $6.05 \times 10^{10}$& 	2.64&	1.00 & 479471 \\
west0132 & 132 &	$4.21\times 10^{11}$&	$2.20 \times 10^{11}$	& 1.12	& $2.20 \times 10^{11}$	&7.49 & 1.00 & 3619199\\
fs\_183\_3	& 183 & $3.27 \times 10^{13}$&	$2.39 \times 10^{13}$	& 1.00	& $2.39 \times 10^{13}$	& $4.73 \times 10^1$	& $6.74 \times 10^{-1}$ & 267\\
    \end{tabular}
    }
    \label{tab:sparsesuit_problems}
\end{table}

\begin{table}[]
\caption{As in Table~\ref{tab:results_uL_single}, but for SuiteSparse problems: $u_L$ is set to single, $u_R$ is set to double. }
    \centering
       \scalebox{0.64}{
    \begin{tabular}{c|c|c|c|c|c|c|c|c}
problem & IC & BE &	FE	& 	$\zeta$&	$\Vert \Zbar_k \Vert \Vert M_R(\xbar_k - \xbar_0) \Vert$	& $\psi_A$&	$\psi_L$	& $\rho$ \\
\hline
rajat14	& 3	& $7.89 \times 10^{-13}$	& $1.84 \times 10^{-6}$	&	$ 1.18 \times 10^{-11}$&	$5.14 \times 10^2$&	$1.42 \times 10^{-2}$&	$1.06 \times 10^{-4}$	&$3.11 \times 10^{-11}$\\
arc130	& 3	& $1.73 \times 10^{-18}$& 	$2.14 \times 10^{-8}$&	$1.41 \times 10^{-16}$&	$ 1.02\times 10^{6}$	&$5.99 \times 10^{-6}$&	$4.80 \times 10^{-10}$	&$4.02 \times 10^{-6}$\\
west0132	& 4	& $2.30 \times 10^{-17}$	& $2.93 \times 10^{-6}$	&	$4.84 \times 10^{-16}$	&$6.46 \times 10^4$	& $1.87 \times 10^{-5}$& 	$6.05 \times 10^{-9}$	&$3.00 \times 10^{-7}$\\
fs\_183\_3	& 3 & 	$2.41 \times 10^{-20}$&	$1.31 \times 10^{-8}$&	$ 1.31 \times 10^{-16}$&	$ 1.53 \times 10^{5}$&	$9.38 \times 10^{-12}$&	$5.32 \times 10^{-13}$	& $9.17 \times 10^{-4}$\\
    \end{tabular}
    }
    \label{tab:sparsesuit_uL_single}
\end{table}

\begin{table}[]
    \caption{As in Table~\ref{tab:results_uR_single}, but for SuiteSparse problems: $u_L$ is set to double, $u_R$ is set to single.}
    \centering
     \scalebox{0.64}{
    \begin{tabular}{c|c|c|c|c|c|c|c|c}
problem & IC & BE &	FE	& 	$\zeta$&	$\Vert \Zbar_k \Vert \Vert M_R(\xbar_k - \xbar_0) \Vert$	& $\psi_A$&	$\psi_L$	& $\rho$ \\
\hline
rajat14	& 3	& $1.25 \times 10^{-19}$	& $6.19 \times 10^{-15}$	&	$1.19 \times 10^{-16}$	& $2.97 \times 10^2$ &	$7.46 \times 10^{-3}$&	$1.18 \times 10^{-4}$	&$1.15 \times 10^{-6}$\\
arc130	& 5&	$3.43 \times 10^{-22}$&	$1.83 \times 10^{-16}$	&	$ 1.11\times 10^{-16}$&	$1.02 \times 10^{6}$&	$2.14 \times 10^{-6}$	&$5.51 \times 10^{-11}$	& $6.49 \times 10^{-4}$\\
west0132 &	5&	$1.48 \times 10^{-21}$	& $4.80 \times 10^{-15}$	&	$ 1.13 \times 10^{-16}$	& $ 6.39\times 10^4$	& $1.98 \times 10^{-5}$&	$6.16 \times 10^{-9}$&	$2.04 \times 10^{-1}$\\
fs\_183\_3	& 3	& $1.52 \times 10^{-27}$&	$ 1.05 \times 10^{-15}$	&	$1.25 \times 10^{-16}$&	$ 1.47 \times 10^{5}$&	$1.87 \times 10^{-11}$&	$1.71 \times 10^{-12}$&	$1.38$ \\
    \end{tabular}
    }
    \label{tab:sparsesuit_uR_single}
\end{table}

\section{Concluding remarks}\label{sec:conclusions}

In light of great community focus on mixed precision computations, we analyzed a variant of split-preconditioned FGMRES that allows using different precisions for computing matrix-vector products with the coefficient matrix $A$ (unit roundoff $u_A$), left-preconditioner $M_L$ (unit roundoff $u_L$), right-preconditioner $M_R$ (unit roundoff $u_R$), and other computations (unit roundoff $u$). A backward error of a required level can be achieved by controlling these precisions. 

Our analysis and numerical experiments show that  
the precision for applying $M_L$ must be chosen in relation to $u$, $u_A$, and the required backward and forward errors, because $u_L$ heavily influences the achievable backward error. 
We can be more flexible when choosing $u_R$ as it does not influence the backward error directly. Our analysis holds under a sufficient but not necessary assumption on $u_R$ in relation to $M_R$. As long as $M_R$ is not singular in precision $u_R$ (note that scaling strategies may be used to ensure this), setting $u_R$ to a low precision is sufficient. Very low precisions $u_L$ and $u_R$ may delay the convergence iteration wise, yet setting $u_L \leq u$ or $u_R \leq u$ does not improve the convergence in general. Note that these conclusions apply to the full left- and right-preconditioning cases as well.

We observe that the forward error is determined by the backward error and the condition number of the left-preconditioned coefficient matrix. This motivates concentrating effort on constructing an appropriate left-preconditioner when aiming for a small forward error: the preconditioner should reduce the condition number sufficiently and needs to be applied in a suitably chosen precision.

\appendix
\section{Proof of Theorem~\ref{th:backward_error}}\label{app:backward_err_proof}
The analysis closely follows \cite{arioli2007note} and \cite{arioli2009using}, and thus we provide the important results for each stage rather than the step-by-step analysis. 

\subsection{Left preconditioner} We start by accounting for the effect of $M_L$. 

\subsubsection{Stage 1: MGS} In this stage, we use precisions $u$, $u_A$ and $u_L$. MGS is applied to 
\begin{equation*}
\Cbar^{(k)} = \begin{bmatrix} fl(M_L^{-1} \rbar_0) &  fl(\Atilde \Zbar_k) \end{bmatrix}.    
\end{equation*}
MGS returns an upper triangular $\Rbar_k$ and there exists an orthonormal $\Vhat_{k+1}$, that is $\Vhat_{k+1}^T \Vhat_{k+1} = I_{k+1}$, such that
\begin{gather}
 \begin{bmatrix} \btilde - \Atilde \xbar_0  &  \Atilde \Zbar_k  \end{bmatrix} \vspace{-1pt} + \vspace{-1pt}  \begin{bmatrix} f_1 + f_2 +f_3 & F_k^{(1)} + F_k^{(2)} \end{bmatrix}\vspace{-1pt}  = \vspace{-1pt} \Vhat_{k+1} \Rbar, \nonumber \\
 \Vert  f_1 \Vert  \leq  (u_A \psi_A + u_L \psi_L) \Vert \Atilde \Vert \Vert  \xbar_0 \Vert , \label{eq:f1_def}\\
\Vert f_2 \Vert \leq  c_1(n) u_L \Vert E_L M_L \Vert \Vert  \btilde \Vert + u \left( \Vert \btilde \Vert + (1+ u_A \psi_A + u_L \psi_L) \Vert \Atilde \Vert \Vert \xbar_0 \Vert \right), \label{eq:f2_bound} \\
 \Vert f_3 \Vert \leq c_2 (n) u \left( \Vert \btilde - \Atilde \xbar_0  \Vert  + (u_A \psi_A + u_L \psi_L) \Vert \Atilde \Vert \Vert \xbar_0 \Vert +  c_1(n) u_L \Vert E_L M_L \Vert \Vert  \btilde \Vert\right), \label{eq:f3_bound} \\
 \Vert F_k^{(1)} \Vert \leq  (u_A \psi_A + u_L \psi_L) \Vert \Atilde \Vert \Vert \Zbar_k \Vert, \label{eq:Fk1_bound}\\
 \Vert F_k^{(2)} \Vert \leq c_3(n,k) u\left( \Vert \Atilde \Zbar_k  \Vert +   (u_A \psi_A + u_L \psi_L) \Vert \Atilde \Vert \Vert \Zbar_k \Vert \right).  \label{eq:Fk2_bound}
\end{gather}
Here $f_1$ is the error in computing the matrix vector product $\Atilde \xbar_0$ and $f_2$ accounts for computing $M_L^{-1}b$ and adding it to the computed $\Atilde \xbar_0$. Error $F_k^{(1)}$ comes from computing $\Atilde \Zbar_k$. $f_3$ and $ F_k^{(2)}$ arise in the MGS process.

\subsubsection{Stage 2: Least squares} The least squares problem is solved using precision $u$. 

From the analysis of \cite{arioli2007note}, under assumptions \eqref{eq:assumptions_nu_ukappaC} and \eqref{eq:assumption_sk} the norm of the residual of the least-squares problem \eqref{eq:LS_problem}
\begin{equation*}
    \alpha_k = \Vert \betabar e_1 + g^{[k]} - (\Hbar_k +\Delta \Hbar_k) \ybar_k  \Vert
\end{equation*}
monotonically converges to zero for a finite $k \leq n$.  
We can express $\alpha_k$ in the following way: 
\begin{equation}
    \alpha_k =  \Vert \btilde - \Atilde \xbar_0 + \delta \rtilde_0 - \Atilde (\Zbar_k + \Zhat_k ) \ybar_k  \Vert,\label{eq:alphak_ZZ}
\end{equation}
where
\begin{gather}
    \delta \rtilde_0 = f_1 + f_2 + f_3  + \Vhat_{k+1} g^{[k]}, \label{eq:delta_r0_def} \\
    \Zhat_k = \Atilde^{-1} \left( F_k^{(1)} +  F_k^{(2)} + \Vhat_{k+1} \Delta \Hbar_k \right), \label{eq:Zhat_k_def}\\
    \Vert g^{[k]} \Vert \leq c_5(k) u \Vert \btilde - \Atilde \xbar_0 \Vert + c_5(k) u (u_A \psi_A + u_L \psi_L) \Vert \Atilde \Vert \Vert \xbar_0 \Vert + c_6(n,k) u u_L \Vert E_L M_L \Vert \Vert  \btilde \Vert, \label{eq:gk_bound_fin}\\
    \Vert \Delta \Hbar_k  \Vert \leq c_4(k) u \Vert \Atilde \Zbar_k \Vert + c_7(n,k) u (u_A \psi_A + u_L \psi_L) \Vert \Atilde \Vert \Vert \Zbar_k \Vert.\label{eq:deltaHk_bound}
\end{gather}

\subsubsection{Stage 3: Computing $\xbar_k$} When certain conditions on the residual norm are satisfied, precision $u$ is used to compute $x_k$ as
\begin{gather}
    \xbar_k = \xbar_0 + \Zbar_k \ybar_k + \delta x_k, \label{eq:xbar_k} \\
    \Vert \delta x_k \Vert \leq c_8(k) u \Vert \Zbar_k \Vert \Vert \ybar_k \Vert + u \Vert \xbar_0 \Vert. \label{eq:bound_deltax_k}
\end{gather}
Using this to eliminate $\Zbar_k \ybar_k$ in \eqref{eq:alphak_ZZ}, then applying the reverse triangle inequality to bound $\Vert \btilde - \Atilde  \xbar_k  \Vert$ and bounding $\Vert \delta \rtilde_0 \Vert$, $\Vert \Atilde \delta x_k \Vert$ and $\Vert \Atilde  \Zhat_k \ybar_k \Vert$ gives
\begin{align}
     \Vert \btilde - \Atilde  \xbar_k  \Vert
 \leq & \, c_{11}(n,k) \Big(\big( u + (1+u)( u_A \psi_A + u_L \psi_L)\big) \Vert \Atilde \Vert \left(  \Vert \xbar_0 \Vert + \Vert \Zbar_k \Vert \Vert \ybar_k \Vert \right) \nonumber\\
 +  & \, \left( u + u_L (1+u)\Vert E_L M_L \Vert\right) \Vert  \btilde \Vert \Big). \label{eq:residual_norm_with_yk}
\end{align}
We eliminate $\Vert \ybar_k \Vert $ from the bound in the following section.

\subsection{Right preconditioner} We now extend the analysis to account for the effect of applying $M_R$. Under assumption \eqref{eq:Mr_assump} $\Zbar_k$ is computed such that
\begin{equation}\label{eq:Z_computed}
    \Zbar_k = M_R^{-1} \Vbar_k + \Delta M_R \Vbar_k, 
\end{equation}
where $\Vert \Delta M_R \Vert \leq c_{12}(n) u_R \Vert E_R \Vert$. 
Then we can obtain
\begin{align}
\Vert \ybar_k \Vert \leq & \, 1.3   \left( \Vert  M_R (\xbar_k - \xbar_0) \Vert +  \Vert M_R \Vert \Vert \delta x_k \Vert + \Vert M_R \Vert \Vert \Delta M_R \Vbar_k \Vert \Vert \ybar_k \Vert \right) \nonumber \\
\leq & \, 1.3 c_{13}(n,k)  \big( \Vert M_R( \xbar_k - \xbar_0) \Vert +  u \Vert M_R \Vert \Vert \Zbar_k \Vert \Vert \ybar_k \Vert + u \Vert M_R \Vert \Vert \xbar_0 \Vert  \label{eq:yk_norm} \\
 + & u_R \Vert M_R \Vert \Vert E_R \Vert \Vert \ybar_k \Vert \big). \nonumber
\end{align}
Under assumption \eqref{eq:rho_def_assump}
\begin{equation*}
    \Vert \ybar_k \Vert \leq \frac{1.3 c_{13}(n,k)}{1 - \rho} \left( \Vert M_R (\xbar_k - \xbar_0) \Vert + u \Vert M_R \Vert \Vert \xbar_0 \Vert \right).
\end{equation*}
Using this in \eqref{eq:residual_norm_with_yk} 
and dropping the terms $u^2$, $u u_L$ and $u u_A$ gives the required result.

\bibliographystyle{siam}
\bibliography{paper}

\end{document}